\documentclass[aos,MSNbibl,nameyear,dvips]{arximspdf}
\usepackage{dcolumn}

\doi{10.1214/11-AOS911}
\volume{39}
\issue{6}
\pubyear{2011}
\firstpage{3092}
\lastpage{3120}

\makeatletter

\def\bsuffix #1{#1}

  \let\sv@tabnotetext\tabnotetext
  \let\sv@tabnotemark@fmt\tabnotemark@fmt
   \long\def\legend#1{{\let\tabnote@indent\leavevmode\sv@tabnotetext[]{}{#1}}}

\newcolumntype{d}[1]{D{.}{.}{#1}}
\newcolumntype{k}[1]{D{,}{}{#1}}

\newtheorem{theorem}{Theorem}[section]
\newtheorem{lemma}{Lemma}[section]

\newproclaim{condition}{Condition}

\newcommand{\bo}{{\mathbf{o}}}
\newcommand{\bu}{\mathbf{u}}
\newcommand{\bz}{\mathbf{z}}
\newcommand{\bA}{\mathbf{A}}
\newcommand{\bB}{\mathbf{B}}
\newcommand{\bE}{\mathbf{E}}

\newcommand{\bQ}{\mathbf{Q}}
\newcommand{\bS}{\mathbf{S}}
\newcommand{\bT}{\mathbf{T}}
\newcommand{\bV}{\mathbf{V}}
\newcommand{\bX}{\mathbf{X}}
\newcommand{\bY}{\mathbf{Y}}
\newcommand{\bfv}{\mathbf{v}}
\newcommand{\bfx}{\mathbf{x}}
\newcommand{\bfe}{\mathbf{e}}
\newcommand{\bfa}{\mathbf{a}}
\newcommand{\bfb}{\mathbf{b}}

\newcommand{\bbeta}{{\bolds\beta}}
\newcommand{\brho}{{\bolds\rho}}
\newcommand{\bxi}{{\bolds\xi}}
\newcommand{\bSigma}{{\bolds\Sigma}}
\newcommand{\var}{\operatorname{var}}

\makeatother

\begin{document}
\begin{frontmatter}

\title{Regularization for Cox's proportional hazards model with
NP-dimensionality\thanksref{T1}}
\runtitle{Regularization for Cox's model}

\thankstext{T1}{Supported
in part by NSF Grants DMS-07-04337, DMS-07-14554 and DMS-09-06482 and
NIH Grant R01-GM072611. Majority of the work was done when
Jelena Bradic was a graduate student at Princeton University.}

\begin{aug}
\author[A]{\fnms{Jelena} \snm{Bradic}\corref{}\ead[label=e1]{jbradic@math.ucsd.edu}},
\author[B]{\fnms{Jianqing} \snm{Fan}\ead[label=e2]{jqfan@princeton.edu}}
\and
\author[C]{\fnms{Jiancheng} \snm{Jiang}\ead[label=e3]{jjiang1@uncc.edu}}
\runauthor{J. Bradic, J. Fan and J. Jiang}
\affiliation{University of California, San Diego,
Princeton University\break and University of North Carolina at Charlotte}
\address[A]{J. Bradic \\
Department of Mathematics \\
University of California, San Diego\\
La Jolla, California 92093-0112 \\
USA\\
\printead{e1}}
\address[B]{J. Fan \\
Department of Operations Research\\
\quad and Financial Engineering \\
Princeton University\\
Bendheim Center for Finance \\
Princeton, New Jersey 08544\\
USA\\
\printead{e2}}
\address[C]{J. Jiang \\
Department of Mathematics and Statistics\\
University of North Carolina at Charlotte\\
Charlotte, North Carolina 28223\\
USA\\
\printead{e3}} 
\end{aug}

\received{\smonth{10} \syear{2010}}
\revised{\smonth{3} \syear{2011}}

%
\begin{abstract}
High throughput genetic sequencing arrays with thousands of
measurements per sample and a great amount of related censored clinical
data have increased demanding need for better measurement specific
model selection. In this paper we establish strong oracle properties of
nonconcave penalized methods for \textit{nonpolynomial} (NP)
dimensional data with censoring in the framework of Cox's proportional
hazards model. A class of folded-concave penalties are employed and
both LASSO and SCAD are discussed specifically. We unveil the question
under which dimensionality and correlation restrictions can an oracle
estimator be constructed and grasped. It is demonstrated that
nonconcave penalties lead to significant reduction of the
``irrepresentable condition'' needed for LASSO model selection
consistency. The large deviation result for martingales, bearing
interests of its own, is developed for characterizing the strong oracle
property. Moreover, the nonconcave regularized estimator, is shown to
achieve asymptotically the information bound of the oracle estimator.
A coordinate-wise algorithm is developed for finding the grid of
solution paths for penalized hazard regression problems, and its
performance is evaluated on simulated and gene association study
examples.
\end{abstract}

%

\end{frontmatter}

\section{Introduction}

A central theme in high-dimensional data analysis is
efficient discovery of sparsity patterns. For such data, where
dimensionality possibly grows exponentially faster than the sample
size, sparsity structures are imposed as means of recovering important signals.
Under the linear regression model framework, various methods ranging
from regularized to marginal regressions and graphical models have been
effectively proposed for identification, reconstruction and estimation
of the unknown sparse regression parameters.

With increasing understanding of sparse recovery in these novel
high-dimen\-sional spaces, more and more attention is paid to efficient
discovery of sparsity patterns for ultra-high dimensional data and
great progress has been made in the least squares setting.
For example,
\citet{MB06}, \citet{ZY06} and \citet{ZH08}
investigated model selection consistency of LASSO when the number of
variables is of a greater order than the sample size and
\citet{CT07} introduced the Dantzig
selector specifically to handle the NP-dimensional variable selection
problem, and \citet{BTW07},
\citet{BRT09}, \citet{vGB09}, \citet{K09}, \citet{MY09}, \citet{MM10}, among
others, showed their asymptotic or finite sample oracle risk properties
for fixed or random ill-posed designs.
Various versions of the ``restricted eigenvalue condition,'' ``sparse
Riesz condition'' or ``incoherence condition'' that exclude high
correlations among variables play a key role here.
On the other hand, when the LASSO estimator does not satisfy some of
these conditions,
it often selects a model which is overly dense in its effort to relax
the penalty on the relevant coefficients [\citet{FL01}, \citet{Z09},
\citet{ZH08}].
Hence, nonconvex penalties [\citet{FL01}] are proposed where \citet{Z09}
pioneered the work with NP-dimensionality and demonstrated its
sign consistency for $p \gg n$ and its advantages over LASSO in the
sense of attaining minimax convergence rates.
\citet{LF09} and \citet{FL10} made important connections between finite
sample and asymptotic oracle properties using folded-concave
penalties for the penalized least
squares estimator with NP-dimensionality.
Although extensive work has been done for linear regression models,
censored survival data have been left greatly unexplored for $p \gg n$.
%

Extending oracle results to censored data with NP-dimensionality
presents a tremendous novel challenge, and, to the best of our
knowledge, there is no previous work on this topic.
The extensions to LASSO and SCAD algorithms for survival data were
successfully proposed by \citet{T97} and \citet{FL02}, respectively, but
both algorithms were theoretically tested only when $p \ll n$.
In recent papers, \citet{J09}, \citet{WNZZ09} and \citet{DML10} addressed
the problem in accelerated failure time models, Cox's model and
semiparametric relative risk models by combining the LASSO, group LASSO
and adaptive LASSO penalties, but, likewise, they only discussed the
case of $p \ll n$.

Motivated by the growing importance of gene selection problems, in
this paper we go one step further and address the problem of existence
of an oracle estimator\vadjust{\goodbreak} and regularization estimator under an ultra-high
dimensionality setting, where the full dimensionality might grow
exponentially or nonpolynomially fast with the sample size, in order
of $\log p =O(n^{\delta})$ for some $\delta>0$,
and the intrinsic dimensionality goes to infinity, in order of
$s=O(n^{\alpha})$ for $\alpha\in(0,1)$. We develop a strong oracle
argument, which shares the spirit of \citet{FL02}, but guarantees that
the folded-concave penalized
partial likelihood estimator is equal to the oracle one, with
probability tending to 1. A similar strong oracle argument was
developed by \citet{KCO08} and \citet{BFW09} in the contexts of linear
regression models. Extending such results to Cox's proportional hazards
model is a new exceptional challenge due to its nature of censoring and
NP-dimensionality.


\subsection{Model setup}

We consider multivariate data $\{(\bX_i, T_i)\}_{i=1}^n$, which form an
i.i.d. sample from the population $(\bX,\bT)$, where $\bX
_i=(X_{i1},\ldots,X_{ip})^T$ is a column vector of covariates for the
$i$th individual. For a variety of reasons
not all survival times $(T_i)_{i=1}^n$ are fully observable. The
independent right censoring scheme is considered where i.i.d. censoring
times $(C_i)_{i=1}^n$ are conditionally independent of survival times
given covariates $\{\bX_i\}_{i=1}^n$. Hence, we work with i.i.d. sample
$
\{(\bX_i,Z_i,\delta_i)\}_{i=1}^{n},
$
where $Z_i=\min(T_i,C_i)$ and $\delta_i= \mathbf{1}\{T_i \leq C_i\}$
are event times and censoring indicator, respectively.

The conditional hazard rate function of $T$ given $\bX=\bfx $
is denoted
by~$\lambda(t|\bfx )$.
Cox's proportional hazards model assumes that
%
\begin{equation}\label{eq1}
\lambda(t|\bX)=\lambda_0(t) \exp(\bbeta^T \bX),
\end{equation}
where the baseline hazard rate $\lambda_0(t)$ is a nuisance function.
Let $t_1<\cdots< t_N$ denote the ordered failure times and $(j)$
denote the label of the item failing at $t_j$. Denote by
$
\mathcal{R}_j=\{i \in\{1,\ldots, n\}\dvtx Z_i \geq t_j\}
$
the risk set at time $t_j$
and by $\Lambda_0(t)=\int_0^t \lambda_0(u)\,du$ the cumulative baseline
hazard function.

Following the approach of nonparametric maximum likelihood estimation,
the ``least informative'' nonparametric modeling of $\Lambda_0(t)$
assumes that $\Lambda_0(t)$ has a jump of size $\theta_j$ at the
failure time $t_j$:
$
\Lambda_0(t;\theta)=\sum_{j=1}^N \theta_j \mathbf{1}\{t_j \leq t\}.
$
If we use the Breslow MLE\vspace*{1pt}
$
\hat{\theta}_j^{-1} =\sum_{i \in\mathcal{R}_j} \exp(\bbeta^T
\mathbf{X}_i)
$,
then the penalized Cox's log partial likelihood becomes [\citet{FL02}]
%
\begin{equation}\label{eq3}
Q_n(\bbeta)- n \sum_{k=1}^p p_{\lambda_n}(|\beta_k|),
\end{equation}
where
$Q_n(\bbeta)=\sum_{j=1}^N
\{ \bbeta^T \mathbf{X}_{(j)} - \log( \sum_{i \in
\mathcal
{R}_j} \exp(\bbeta^T \mathbf{X}_i ))
\}$,
$p_{\lambda_n}(\cdot)$ is a penalty function, and $\lambda_n$ is a
nonnegative regularization parameter.
Note that the covariate vector $\bX$ may be time dependent and
incorporated in the standard way
in model (\ref{eq1}) through
\[
\lambda(t|\bX(t))=\lim_{\Delta t \to0}P\{t\le T\le t+\Delta t|T\ge t,
\bX(t)\}/ \Delta t=\lambda_0(t) \exp(\bbeta^T \bX(t) ).
\]
Note that from hereon we will be working with the time-dependent left
continuous covariate vector $\bX(t)$.

\subsection{Counting process representation}


Let $N_i(t)=1\{Z_i \leq t, \delta_i=1 \}$,\break $ \bar{N}(t) = \sum _{i=1}^n
N_i(t)$ and $Y_i(t)=1\{Z_i \geq t\}$. Note that the process $\bY
(t)=(Y_1(t),\ldots,\allowbreak Y_n(t))^T$ is assumed to be left continuous with
right-hand limits and satisfies $P(\bY(t)=1, 0\le t\le\tau)>0$. Using
the counting process notation, one can rewrite the log partial
likelihood $Q_n(\bbeta)$ for model (\ref{eq3}) as
\[
Q_n(\bbeta)=\sum_{i=1}^n \int_{0}^{\tau} \bigl\{
\bbeta^T \bX_i(t) -\log\bigl( S^{(0)}_{n}(\bbeta,t) \bigr)\bigr\}
\,d {N}_i(t) ,
\]
where and hereafter
$\tau$ is the study ending time,
and
\[
S^{(\ell)}_{n}(\bbeta,t)=n^{-1}\sum_{i=1}^n Y_i(t)\{\bX_i(t)\}
^{\otimes
\ell}\exp( {\bbeta}^T \mathbf{X}_i(t)),\qquad
\ell=0,1,2,
\]
%
with $\otimes$ denoting the outer product.
Thus, the penalized log partial likelihood becomes
%
\begin{equation}\label{eq4}\quad
\mathcal{C}(\bbeta,\tau) \equiv\sum_{i=1}^n \int_{0}^{\tau} \bigl\{
\bbeta^T \bX_i(t) -\log\bigl( S^{(0)}_{n}(\bbeta,t) \bigr)\bigr\}\, d {N}_i(t)
- n \sum_{j=1}^p p_{\lambda_n}(|\beta_j|).
\end{equation}
Define the sparse estimator $\hat{\bbeta}$ as the maximizer of
$\mathcal{C}(\bbeta,\tau)$ over $\bbeta\in\Omega_p$,
where~$\Omega_p$ is the parameter space which is a compact subset of
$R^p$ and contains the true value of $\bbeta$.
Note that $N_i(t)$ is a counting process with intensity process
$\lambda_i(t,\bbeta)=\lambda_0(t)Y_i(t) \exp\{\bbeta^T X_i(t)\}$, which
does not admit jumps at the same time as $N_j(t)$ for $j \neq i$.
Denote by $\bbeta^*$ the true value of $\bbeta$
and
$\Lambda_i(t)=\int_0^t \lambda_i(u,\bbeta^*)\,du$.
Then $M_i(t)=N_i(t)-\Lambda_i(t)$ is an orthogonal local square
integrable martingale
with respect to filtration
\[
{\mathcal F}_{t,i}=\sigma\{N_i(u), \bX_i(u^+), Y_i(u^+), 0\le u\le t\},
\]
that is, $\langle M_i(t),M_j(t) \rangle=0$ for $i\neq j$.
Let ${\mathcal F}_t=\bigcup_{i=1}^n {\mathcal F}_{t,i}$ be the smallest
$\sigma$-algebra containing ${\mathcal F}_{t,i}$.
Then $\bar{M}(t)=\sum_{i=1}^n M_i(t)$ is a martingale with respect to
${\mathcal F}_t$.

\subsection{Choice of the penalty function}

There are many commonly used penalties in the literature, for example,
the $L_2$ penalty used in ridge regression;
the nonnegative garrote as a shrinkage estimation [\citet{YL07}];
the $L_0$ penalty for the best subset selection;
the $L_1$ penalty LASSO [\citet{T96}] as a convex relaxation of the $L
_0$ penalty; the SCAD penalty [\citet{FL01}],
defined via its derivative
$p_{\lambda}'(t) = \lambda\{
I (t \le\lambda) +
\frac{(a\lambda-t)+}{(a -1) \lambda}
I (t > \lambda)
\}
, t \ge0$, for some $a > 2$,
as a folded-concave relaxation of~$L_0$ penalty;
the MCP [\citet{Z09}] penalty.
Recently, a class of penalties bridging $L_0$ and $L_1$ penalties was
introduced in \citet{LF09}. All of these penalties are folded concave
penalties, as noted in \citet{FL01} and \citet{FL10}. As a collection of
nonconvex relaxations of the $L_0$ penalty, they serve as a tool of
allowing bigger correlations among covariates (see Condition \ref
{cond4}) and hence relax significantly the standard ``incoherence
condition'' and control the tail bias of the resulting penalized
estimator (see Theorem~\ref{theo41}).
For any penalty function $p_{\lambda_n}(\cdot)$, let $\rho(t;
\lambda
_n) = \lambda_n^{-1}p_{\lambda_n}(t)$
and write $\rho(t; \lambda_n)$ as $\rho(t)$ for simplicity when there
is no confusion.
According to \citet{FL10}, the folded concave penalties are defined
through the following Condition \ref{cond1}.
%
\begin{condition}\label{cond1}
$ \rho(t; \lambda_n)$ is increasing and concave in $t \in[0,\infty)$
and has a continuous derivative $ \rho'(t; \lambda_n)$
with $ \rho'(0+; \lambda_n)> 0$. In addition, $ \rho'(t; \lambda_n)$
is increasing in $\lambda_n \in(0,\infty)$
and $ \rho'(0+; \lambda_n)$ is
independent of $\lambda_n$.
\end{condition}

Note that most commonly used nonconvex penalties, including SCAD and
MCP ($a \ge1$), satisfy Condition \ref{cond1}.
We will employ the folded concave penalties to increase flexibility of our
method.
LASSO penalty as a convex function falls at the boundary of penalties
in Condition \ref{cond1}, and our results will be applicable for LASSO
penalty as well.

The rest of the paper is organized as follows. In Section \ref{sec2}
we deal with identification problem of the penalized estimator $\hat
{\bbeta}$, which is key to the proof of oracle results. A compelling
large deviation result is derived for divergence of a martingale from
its compensator in Section \ref{sec3}. In Section \ref{sec4} we work
out the new strong oracle property and its implications for LASSO and
SCAD and asymptotic properties of the proposed estimator. In Section~\ref{sec7} we propose an iterative coordinate ascent algorithm (ICA)
and examine a~thorough simulation example; see Section~\ref
{subsec72}. The gene association study is done in Section \ref
{subsec73} where the non-Hodgkin's lymphoma dataset of \citet{D04} is
analyzed. Technical lemmas and proofs are collected in the Appendix and
in the supplementary material [\citet{BFJ11}].

\section{Identification} \label{sec2}

This section gives the appropriate necessary and sufficient conditions
on the existence of estimator $\hat{\bbeta}$. We can always assume that
the true parameter $\bbeta^*$ can be arranged as
$\bbeta^*=(\bbeta_1^{*T},\mathbf{0}^T)^T$,
with $\bbeta_1^*\in\Omega_s$ being a vector of nonvanishing elements
of $\bbeta^*$,
where $\Omega_s=\Omega_p\cap R^s$.

Throughout the paper
the following notation on a vector/matrix norm is used.
Denote by
$\lambda_{\min}(\bB)$ and $\lambda_{\max}(\bB)$
the minimum and maximum eigenvalues of a symmetric matrix $\bB$, respectively.
We also use $\lambda(\bB)$ to denote any eigenvalue of~$\bB$.
Let \mbox{$\|\cdot\|_q$} be the $L_q$ norm of a vector or matrix.
Then
for a $s\times s$ matrix $\bA$,
$\|\bA\|_{\infty}=\max\{\sum_{k=1}^s | (\bA)_{jk}|\dvtx 1\le j\le s\}$
and $\|\bA\|_2=\{\lambda_{\max}(\bA^T\bA)\}^{1/2}$.
We also let
$\sigma(\bA)$ be the set consisting of all of eigenvalues of $\bA$,
and let $r_{\sigma}(\bA)=\max\{|\lambda|\dvtx \lambda\in\sigma(\bA
)\}$
be the spectral radius of $\bA$.
If $\bA$ is symmetric, then $r_{\sigma}(\bA)=\|\bA\|_2$.

Since no concavity is assumed for the penalized log partial likelihood
(\ref{eq4}),
it is difficult, in general, to study the global maximizer of the
penalized likelihood.
One useful index controlling the convexity of the whole optimization
problem (\ref{eq4}) is the following ``local concavity'' of the penalty
function $\rho(\cdot)$ at $\bfv = (v_1,\ldots, v_s)^T \in
R^s$ with $\|
\bfv \|_0 = s $,
%
\begin{equation}
\kappa(\rho, \bfv )=\lim_{\varepsilon\to0+} \max_{1\leq
j\leq s}\sup
_{t_1< t_2 \in(|v_j|-\varepsilon, |v_j|+\varepsilon)}
-\frac{\rho'(t_2)-\rho'(t_1)}{t_2-t_1},
\end{equation}
which is defined in \citet{LF09} and
shares similar spirit to the ``maximum concavity'' of $\rho$ in \citet{Z09}.
Since $\rho$ is concave on $(0,\infty)$, $\kappa(\rho, \bfv ) \ge0$.
For LASSO penalty $\kappa(\rho,\bfv )=0$,
whereas for the SCAD penalty
\[
\kappa(\rho, \bfv )=\cases{
(a - 1)^{-1}\lambda^{-1}, &\quad if there exists a
$v_j$ such that $\lambda\le|v_j|\le a\lambda$;\cr
0, &\quad otherwise.}
\]

Let ${\bbeta}_1$ be a subvector of ${\bbeta}$ formed by all nonzero
components and $s=\dim({\bbeta}_1)$. Denote by $\bS_{i}$ the subvector
of $\bX_i$ with same indexes as $\hat {\bbeta}_1$ in~$\hat{\bbeta}$ and
by $\bQ_{i}$ the complement to $\bS_{i}$. For $\bfv
=(v_1,\ldots,v_s)^T\in R^s$, let $\brho'(\bfv )
=(\rho'(v_1),\ldots,\rho'(v_s))^T$ and $\operatorname{sgn}(\bfv )
=(\operatorname{sgn}(v_1),\ldots,\operatorname{sgn}(v_s))^T$.
Partition
$S_n^{(1)}(\bbeta,\break t)=[{S_{n1}^{(1)}(\bbeta,t)},{S_{n2}^{(1)}(\bbeta,t)}]$
and
\[
S_n^{(2)}(\bbeta,t)=\pmatrix{
S_{n11}^{(2)}(\bbeta,t) & S_{n12}^{(2)}(\bbeta,t)\vspace*{2pt}\cr
S_{n21}^{(2)}(\bbeta,t) & S_{n22}^{(2)}(\bbeta,t)}
\]
according to the partition of $\bbeta=(\bbeta_1^T,\bbeta_2^T)^T$, so
that $S_{n1}^{(1)}(\bbeta,t)$ is a \mbox{$\dim(\bbeta_1)\times1$} vector and
$S_{n11}^{(2)}(\bbeta,t)$ is a $\dim(\bbeta_1)\times\dim(\bbeta _1)$
matrix. Let
$\bE_n^{(1)}(\bbeta,t)=S_{n1}^{(1)}(\bbeta,\allowbreak t)/ S_n^{(0)}(\bbeta,t)$,
$\bE_n^{(2)}(\bbeta,t)=S_{n2}^{(1)}(\bbeta,t)/\allowbreak S_n^{(0)}(\bbeta,t)$,
$\bE_n(\bbeta,t)=S_{n}^{(1)}(\bbeta,t)/S_n^{(0)}(\bbeta,t)$,
$\bV(\bbeta,t) = {S_{n11}^{(2)}(\bbeta,t)}/{S_n^{(0)}(\bbeta,t)}
-({S_{n1}^{(1)}(\bbeta,t)}/{S_n^{(0)}(\bbeta,t) } )^{\otimes2}$ and
$\bV(\bbeta_1,t)=\break\bV((\bbeta_1,\mathbf{0}),t)$.

The following theorem provides a sufficient condition on the strict
local maximizer of ${\mathcal C}(\bbeta,\tau)$. Proof is relegated to
the supplementary material [\citet{BFJ11}].
%
\begin{theorem}\label{theo21}
If Condition \ref{cond1} is satisfied, then an estimate $\hat
{\bolds\beta} \in R^p$ is a~strict local maximizer of the
nonconcave penalized log partial likelihood (\ref{eq4}) if
%
\begin{eqnarray}\label{eq5}
&\displaystyle \sum_{i=1}^n \int_{0}^{\tau} \bigl( {\bS}_{i}(t) - \bE
_n^{(1)}(\hat
\bbeta,t) \bigr)\,dN_i(t) - n \lambda_n {\brho}'(|\hat{\bolds
\beta
}_1|) \circ\operatorname{sgn}(\hat{\bolds\beta}_1)= \mathbf{0},
&
%
\\
\label{eq6}
&\displaystyle \|\bz(\hat{\bbeta})\|_{\infty}\equiv\Biggl\| \sum_{i=1}^n \int
_{0}^{\tau
} \bigl({\bQ}_{i}(t) - \bE_n^{(2)}(\hat{\bbeta},t) \bigr)\, d
N_i(t) \Biggr\|
_{\infty} < n \lambda_n \rho'(0+),&
\\
%
\label{eq7}
&\displaystyle \lambda_{\min} \biggl\{ \int_{0}^{\tau}\bV(\hat\bbeta,t) \,d \bar{N}(t)
\biggr\} >
n \lambda_n \kappa(\rho, \hat{\bolds\beta}_1),&
\end{eqnarray}
where $\circ$ is the Hadamard product.
Conversely, if $\hat\bbeta$ is a local maximizer of $\mathcal
{C}(\bbeta
,\tau)$,
then it must satisfy (\ref{eq5})--(\ref{eq7}) with strict inequalities
replaced by nonstrict inequalities.
\end{theorem}

When the LASSO penalty is used, $\kappa(\rho,\bfv )=0$, hence the
condition of nonsingularity for the matrix in (\ref{eq7}) is
automatically satisfied with a nonstrict inequality.
For the SCAD penalty, $\kappa(\rho, \hat{\bolds\beta}_1) = 0$;
that is, (\ref{eq7}) holds with nonstrict inequality,
unless there are some $j$
such that
$\lambda_n<|\hat\beta_j|< a\lambda_n$,
which usually has a small chance.
In the latter case,
$\kappa(\rho, \hat{\bolds\beta}_1)= (a- 1)^{-1}\lambda_n^{-1}$,
and the condition in (\ref{eq7}) reduces to
\[
\lambda_{\min}\biggl\{\int_{0}^{\tau}n^{-1}\bV(\hat\bbeta,t)
\,d\bar
{N}(t)\biggr\}>1/(a-1).
\]
This will hold if a large $a$ is used, due to nonsingularity of the matrix.

It is natural to ask if the penalized nonconcave Cox's log partial
likelihood has a global maximizer. Since $p \gg n$, it is\vspace*{1pt}
hard to show the global optimality of a local maximizer. Theorem
\ref{theo411} in Section \ref{sec4} suggests a condition for
$\hat{\bbeta}$ to be\vspace*{1pt} unique and global. Once the unique maximizer is
available, it will be equal to the oracle one with probability tending
to one exponentially fast, when the effective dimensionality $s$ is
bounded by $O(n^{\alpha})$ for $\alpha<1$ (see Theorem \ref{theo42}).
In this way Theorems \ref{theo21}, \ref{theo411} and \ref{theo42}
address uniqueness of the solution and provide methods for finding the
global maximizer among potentially many. Methodological innovations
among others consist of using equations (\ref{eq5}) and (\ref{eq6}) as
an identification tool to surpass the absence of analytical form of an
estimator $\hat \bbeta$.

%

\section{A large deviation result} \label{sec3}


In view of (\ref{eq5}) and (\ref{eq6}), to study a nonconcave
penalized Cox's partial likelihood estimator $\hat\bbeta$, we need to
analyze the deviation of $p$-dimensional counting process $\int_0^t \{
\bX_i(u)-\bE_n(\bbeta^*,u)\} \,d N_i(u)$ from its compensator $\bA
_i=\int_0^t \{\bX_i(u)-\bE_n(\bbeta^*,u)\} \,d\Lambda_i(u)$. In other
words, we need to simultaneously analyze the deviation of marginal
score vectors from their compensators.
Some conditions are needed for this purpose.
%
\begin{condition}\label{cond3}
There exists a compact neighborhood ${\mathcal B}$ of $\bbeta^*$ that
satisfies each of the following conditions:

\begin{longlist}
\item
There exist scalar, vector and matrix functions $s^{(j)}$ defined on
${\mathcal B}\times[0,\tau]$ such that,
in probability as $n\to\infty$ for $j=0,1,2$,
${\sup_{t\in[0,\tau],{\beta_1}\in{\mathcal B_1}}}\|S_n^{(j)}(\bbeta
_1,\allowbreak t)-s^{(j)}(\bbeta_1, t)\|_{2}\to0$,
for $\mathcal B_1 \in\mathcal{R}^s, \mathcal B_1 \subset\mathcal B$.

\item
The functions $s^{(j)}$ are bounded and $s^{(0)}$ is bounded away from
$0$ on ${\mathcal B}\times[0,\tau]$; for $j=0,1,2$, the family of
functions $s^{(j)}(\cdot, t)$, $0\le t\le\tau$, is an equicontinuous
family at $\bbeta^*$.

\item Let
$\bfe (\bbeta,t)=s^{(1)}(\bbeta,t)/s^{(0)}(\bbeta,t)$,
$\bfv (\bbeta,t)=s^{(2)}(\bbeta,t)/s^{(0)}(\bbeta,t)-\{\bfe(\bbeta,\allowbreak
t)\}^{\otimes2}$
and
$\bSigma_{\beta}(t)=\int_0^t\bfv (\bbeta,u)s^{(0)}(\bbeta
^*,u) \,d\Lambda_0(u)$.
Define
$\bfv (\bbeta_1,t)$
in the same way as for
$\bV(\bbeta_1,t)$ but with $S_n^{(\ell)}$ replaced by $s^{(\ell)}$.
Let
\[
\bSigma_{\beta_1}(t)=\int_0^t\bfv (\bbeta
_1,u)s^{(0)}(\bbeta_1^*,u)\,
d\Lambda_0(u)
\]
%
and $\bSigma_{\beta_1}=\bSigma_{\beta_1}(\tau)$.
Assume that the $s\times s$ matrix
$\bSigma_{\beta^*_1}$ is positive definite for all $n$
and
$\Lambda_0(\tau) < \infty$.


\item Let
$c_n={\sup_{t\in[0,\tau]}}\|\bE_n(\bbeta^*,t)-\bfe (\bbeta
^*,t)\|_{\infty}$
and
$d_n={\sup_{t\in[0,\tau]}}|S_{n}^{(0)}(\bbeta^*,\allowbreak
t)-s^{(0)}(\bbeta^*,t)|$.
The random sequences $c_n$ and $d_n$ are bounded almost surely.
\end{longlist}
\end{condition}

The above conditions, (i)--(iii), agree with the conditions in Section
8.2 of \citet{FH05} for fixed $p$ and in \citet{CFLZ05} for
diverging $p$. Condition (iii) is restricted to hold on the $s$ instead
of usually assumed $p$-dimensional subspace. This is a counterpart of
the similar conditions imposed on the covariance matrix $\bX$ in the
linear regression models [see, e.g., \citet{BTW07},
\citet{vGB09}, \citet{Z09}]. Nonsingularity of the matrix
$\bSigma _{\beta ^*_1}$ in (iii) could have been relaxed toward
restricted eigenvalue properties like those for linear models
[\citet{BRT09}, \citet{K09}] but for easier composure we impose a bit
stronger condition.
Condition (iv) is used to ensure that
the score vector of the log partial likelihood, which is a martingale,
has bounded jumps
and quadratic variation.
By following the discussion on pages 305 and 306 of \citet{FH05}, this
condition is not stringent for i.i.d. samples.

The following Condition \ref{cond2} is coming as a consequence of
martingale representation of the score function for the Cox model, and
it is valuable in analyzing large deviations of counting processes.
%
\begin{condition}\label{cond2}
Let
$\varepsilon_{ij}=\int_0^{\tau}(
X_{ij}(t)-e_j(\bbeta^*,t)
)\, dM_i(t)$,
where
$e_j(\bbeta^*,t)$ is the $j$th component of $\bfe (\bbeta^*,t)$.
Suppose the Cram\'er condition holds for $\varepsilon_{ij}$, that is,
\[
E|\varepsilon_{ij}|^m\le m!M^{m-2}\sigma_j^2/2
\]
for all $j$, where $M$ is a positive constant,
$m\ge2$
and
$\sigma_j^2=\var(\varepsilon_{ij})
<\infty$.
\end{condition}

In linear regression models, the large deviation is established upon
the Cram\'er condition for the covariates. Condition \ref{cond2} takes
a similar role here and can be regarded as an extension to the
classical Cram\'er condition.
Moreover, it is trivially fulfilled if the covariates are bounded. In
that sense it represents a relaxation of typical assumption of bounded
covariates.
Since~$\varepsilon_{ij}$ is a mean zero martingale,
it can be shown that
$\sigma_j^2=E(\varepsilon_{ij}^2)=(\bSigma_{\beta^*})_{jj}$ is the
$j$th diagonal entry of $\bSigma_{\beta^*}$.
Define $\bxi=(\xi_1,\ldots,\xi_p)^T$ to be the score vector of the log
partial likelihood function $Q_n(\bbeta)$,
\[
\bxi=\sum_{i=1}^n \int_{0}^{\tau} \bigl(\bX_{i}(t) -\bE_n(\bbeta
^*,t)\bigr)\, d{N}_i(t).
\]
Since $M_i(t)=N_i(t) -\Lambda_i(t)$ is a martingale with compensator
$\Lambda_i(t)=\break\int_0^t \lambda_i(u,\bbeta^*) \,du$,
we can rewrite $\xi_j$ as
$
\sum_{i=1}^n \int_{0}^{\tau} \{ {X}_{ij}(t)-E_{nj}(\bbeta^*,t)
\} (d M_i(t) + d \Lambda_i(t)),
$
where
$E_{nj}(\bbeta^*,t)$ is the $j$th component of $\bE_{n}(\bbeta^*,t)$.
Note that
$
\sum_{i=1}^n \int_0^{\tau} \{ {X}_{ij}(t)-E_{nj}(\bbeta^*,t)
\} \,d\Lambda_i(t) =0,
$
leading to the representation of the form
\[
\xi_j = \sum_{i=1}^n \int_{0}^{\tau} \{
{X}_{ij}(t)-E_{nj}(\bbeta
^*,t) \} \,d M_i(t).
\]
The following theorem characterizes the uniform deviation of the score
vector~$\bxi$ and is critical in obtaining strong oracle property; see
Theorem \ref{theo42} in Section \ref{sec4}. To the best of our
knowledge there is no similar result in the literature.
%
\begin{theorem}\label{theo31}
Under Conditions \ref{cond3} and \ref{cond2},
for any positive sequence~$\{u_n\}$ bounded away from zero
there exist positive constants $c_0$ and $c_1$ such that
%
\begin{equation}\label{eq31}
P \bigl( |\xi_j|>\sqrt{n}u_n \bigr)
\leq c_0 \exp(- c_1u_n)
\end{equation}
uniformly over $j$, if $v_n=\max_j\sigma_j^2/u_n$ is bounded.
\end{theorem}
\begin{pf}
$\!\!\!$Denote by
$E_{nj}(\bbeta^*,t)$ and $e_j(\bbeta^*,t)$
the $j$th components of $\bE_n(\bbeta^*, t)$ and $\bfe (\bbeta
^*,t)$, respectively.
Then
$\xi_j$ can be written as
\begin{eqnarray*}
\xi_j&=&\sum_{i=1}^n\int_0^{\tau}\bigl(
X_{ij}(t)-e_j(\bbeta^*,t)
\bigr) \,dM_i(t)\\
&&{}-\sum_{i=1}^n\int_0^{\tau}\bigl(
E_{nj}(\bbeta^*,t)-e_j(\bbeta^*,t)
\bigr)\, dM_i(t)\\
&\equiv& \xi_{j1}(\tau)-\xi_{j2}(\tau).
\end{eqnarray*}
To establish the exponential inequality about $\xi_j$,
in the following we will establish the exponential inequalities about
$\xi_{j1}(\tau)$ and $\xi_{j2}(\tau)$.

Note that $\xi_{j1}(\tau)=\sum_{i=1}^n\varepsilon_{ij}$,
where
$\{\varepsilon_{ij}\}_{i=1}^n$ is a sequence of i.i.d. random variables
with mean zero and satisfying\vadjust{\goodbreak}
Condition \ref{cond2}. It follows from the Bernstein exponential
inequality that
%
\begin{equation} \label{eqj1}
P(|\xi_{j1}|>a)\le2 \exp\{-a^2/2(n\sigma_j^2+Ma)\}.
\end{equation}

Note that $\bar{M}(t)$ is a martingale with respect to
${\mathcal F}_t$; it follows that
$\xi_{j2}(t)$
is also a martingale with respect to ${\mathcal F}_t$.
Let
$\bar{N}(t)=\sum_{i=1}^n N_i(t)$. Then
$\Delta\bar{N}(t)=\sum_{i=1}^n \Delta N_i(t)$,
where and thereafter
$\Delta N_i(t)=N_i(t)-N_i(t^-)$ denotes the jump of $N_i(\cdot)$ at
time $t$.
Since no two counting processes $N_i$ jump at the same time, we have
$|\Delta\bar{N}(t)|\le1$.
Let
$\bar{\Lambda}(t)=\sum_{i=1}^n\Lambda_i(t)$.\vspace*{1pt}
By continuity of the compensator $\Lambda_i(t)=\int_0^t\lambda
_i(u,\bbeta^*) \,du$,
$|\Delta\bar{\Lambda}(t)|=0$.
Since $\bar{M}(t)=\bar{N}(t)-\bar{\Lambda}(t)$,
$|\Delta\bar{M}(t)|=|\Delta\bar{N}(t)|\le1$.
Note that $\bY(t)$ and $\bX(t)$ are left continuous in $t$.
It is easy to see that
\begin{eqnarray*}
|\Delta(n^{-1/2}\xi_{j2}(t))|&=& n^{-1/2}|E_{nj}(\bbeta
^*,t)-e_j(\bbeta
^*,t)| \\
&\le& {n^{-1/2}\sup_{t\in[0,\tau]}}\|\bE_{n}(\bbeta^*,t)-\bfe (\bbeta^*,t)\|
_{\infty}\\
&\equiv& n^{-1/2}c_n,
\end{eqnarray*}
which is bounded almost surely by Condition \ref{cond2}(vi).
Note that the predictable quadratic variation of $n^{-1/2}\xi_{j2}(t)$,
denoted by\vadjust{\goodbreak} $\langle n^{-1/2}\xi_{j2}(t)\rangle$, is bilinear and
satisfies that
\begin{eqnarray*}
\langle n^{-1/2}\xi_{j2}(t)\rangle
&=& n^{-1}\int_0^t
\bigl(
E_{nj}(\bbeta^*,u)-e_j(\bbeta^*,u)
\bigr)^2 \,d\langle\bar{M}(u)\rangle\\
&=&\int_0^t\{E_{nj}(\bbeta^*,u)-e_j(\bbeta^*,u)\}^2S_n^{(0)}(\bbeta
^*,u) \,d\Lambda_0(u)\\
&\le& \int_0^t\|\bE_{n}(\bbeta^*,u)-\bfe (\bbeta^*,u)\|
_{\infty
}^2S_n^{(0)}(\bbeta^*,u) \,d\Lambda_0(u)
\equiv b_n^2(t).
\end{eqnarray*}
%
Obviously,
$b_n^2(t)\le b_n^2(\tau)\leq c_n^2 \int_0^{\tau}S_n^{(0)}(\bbeta
^*,t)\,
d\Lambda_0(t)$.
Note that
\[
\int_0^{\tau}S_n^{(0)}(\bbeta^*,t) \,d\Lambda_0(t)
\le\int_0^{\tau}s^{(0)}(\bbeta^*,t) \,d\Lambda_0(t)
+d_n\Lambda_0(\tau).
\]
By Condition \ref{cond3}(ii), (iii) and (vi),
there exist constants $0\le K<\infty$ and $0<b<\infty$, independent
of $j$, such that
$|\Delta(n^{-1/2}\xi_{j2}(t))|\le K$
and
$ \langle n^{-1/2}\xi_{j2}(t)\rangle\le b^2$.
It follows from the exponential inequality for martingales with
bounded jumps
[see Lemma 2.1 of \citet{vG95}] that,
for $u_n>0$,
%
\[
P\bigl\{|\xi_{j2}(\tau)|>\sqrt{n}u_n\bigr\}
= P\{|n^{-1/2}\xi_{j2}(\tau)|>u_n\} \le2\exp\biggl\{-\frac
{u_n^2}{2(Ku_n+b^2)}\biggr\}.
\]
Therefore, by Condition \ref{cond2}(iv),
there exists a constant $c>0$ such that
%
\begin{equation}\label{eqj2}
P\bigl\{|\xi_{j2}(\tau)|>\sqrt{n}u_n\bigr\}
\le2\exp\{-cu_n\}
\end{equation}
uniformly over $j$.
Note that
\[
P\bigl\{|\xi_{j}(\tau)|>\sqrt{n}u_n\bigr\}
\le P\bigl\{|\xi_{j1}(\tau)|>0.5\sqrt{n}u_n\bigr\}+
P\bigl\{|\xi_{j2}(\tau)|>0.5\sqrt{n}u_n\bigr\}.
\]
It follows from (\ref{eqj1}) and (\ref{eqj2}) that
$P\{|\xi_{j}(\tau)|>\sqrt{n}u_n\}$ is bounded by
%
\begin{equation}\label{eqc3}
2\exp\biggl\{-\frac{u_n}{4(2 \sigma
_j^2u_n^{-1}+Mn^{-1/2})}\biggr\}
+2\exp(-0.5cu_n).
\end{equation}
Then there exist positive constants $c_0$ and $c_1$ such that
$P\{|\xi_{j}(\tau)|>\sqrt{n}u_n\}<c_0\exp(-c_1u_n)$
uniformly over $j$, if
$\max_j\sigma_j^2=O(u_n)$.
\end{pf}

Theorem \ref{theo31} represents a uniform, nonasymptotic exponential
inequality for martingales. Compared with other exponential
inequalities [\citet{dlP99}, \citet{JN10},
\citet{vG95}], it is uniform over all components $j$.
Moreover, its independence of dimensionality $p$ proves to be
invaluable for NP variable selection.


\section{Strong oracle property}\label{sec4}

In this section we will prove a strong oracle property result, that is,
that $\hat\bbeta$ is an oracle estimator with overwhelming
probability, and not that it behaves like an oracle estimator [\citet
{FL02}]. We assume that the effective and full dimensionality satisfy
$s=O(n^{\alpha})$ and $\log p = O(n^{\delta})$, for some $\alpha\in
(0,1)$ and $\delta>0$, respectively. This notion of strong oracle
property requires a definition of biased oracle estimator as it was
defined in \citet{BFW09} for the linear regression problem.

Let us define the biased oracle estimator $\hat{\bolds\beta
}{}^{\bo
}=(\hat{\bolds\beta}{}^{\bo T}_1,\mathbf{0}^T)^T$ where $\hat
{\bolds\beta}{}^{\bo}_1$ is a~solution to the $s$ dimensional sub-problem
\[
\mathop{\arg\max}_{\beta_1 \in\Omega_s} \sum_{i=1}^n \int_{0}^{\tau}
\bigl[\bbeta_1^T \bS_i(t) - \log\bigl( S^{(0)}_{n}((\bbeta_1,\mathbf{0}),t) \bigr)
\bigr] \,d N_i(t) - n\lambda_n \sum_{j=1}^s \rho(|\beta_j|;\lambda_n).
\]
That is,
$\hat{\bolds\beta}{}^{\bo}_1 = \arg\max\{\mathcal{C}(\bbeta
_1,\tau
)\dvtx \bbeta_1 \in\Omega_s\}$ with $\mathcal{C}(\bbeta_1,\tau
)=\mathcal
{C}((\bbeta_1,\mathbf{0}),\tau)$.
The estimator
$\hat{\bolds\beta}{}^{\bo}$ is called the biased oracle estimator,
since the oracle knows the true submodel
${\mathcal M}_*=\{j\dvtx \beta_j^*\ne0\}$,
but still applies a penalized method to estimate the nonvanishing coefficients.
%
\begin{theorem}[(Global optimality)]
\label{theo411}
Suppose that
$
\min_{\mathbf{\beta}_1 \in\Omega_s}\lambda_{\min} \{ \int
_{0}^{\tau}\bV(\bbeta_1,\allowbreak t) \,d\bar{N}(t) \} >
n \lambda_n \kappa(\rho, {\bolds\beta}_1)
$
holds almost surely. Then $\hat{\bolds\beta}{}^{\bo}_1$ is a unique
global maximizer of the penalized log-likelihood $\mathcal{C}(\bbeta
_1,\tau)$ in $\Omega_s$.
\end{theorem}

The above theorem could be relaxed to a minimum over the level sets of
Cox's partial likelihood in a similar manner
to Proposition 1 of \citet{FL10}. Its proof is left for the
supplementary material [\citet{BFJ11}].
For LASSO penalty, $\hat{\bolds\beta}{}^{\bo}_1$ is unique and is
the global maximizer, since $\mathcal{C}(\bbeta_1,\tau)$
is strictly concave. In general,
global maximizers are available for SCAD and MCP penalties, if one uses
a large parameter $a$.
In this setting, the biased oracle estimator is unique as a solution
to strictly concave optimization problem.
Note that it still depends on the penalty function.
The biased oracle estimator, by its definition, satisfies only
equation (\ref{eq5}) in Theorem \ref{theo21}.
Since the vanishing component does not need any penalty, the smaller
the penalty the less the bias. In this sense, the biased oracle
estimator with the SCAD penalty has a better performance than the
biased oracle estimator with the LASSO penalty.
The former is asymptotically unbiased [the second term in (\ref{eq5})
is zero], while the latter is not (see Theorems \ref{cor2} and
\ref{cor3}).

In order to establish asymptotic properties of $\hat{\bbeta}_1$, we
need to govern the conditioning number of the $s \times s$ information
matrix $\bSigma_{\beta_1^*}$ through its eigenvalues.
This is done in the following condition.
%
\begin{condition}\label{cond6}
$r_{\sigma}({\bSigma}_{\beta_1^*})=O(1)$
and
$r_{\sigma}({\bSigma}_{\beta_1^*}^{-1})=O(1)$.
\end{condition}

Concerning Condition \ref{cond3}, positive definiteness of ${\bSigma
}_{\beta_1^*}$ is not enough and further bound on its spectrum is
needed. Condition \ref{cond6} is in the same spirit as the partial
Riesz condition
and is weaker than Condition A3 of \citet{CFLZ05},
where Condition \ref{cond6} holds for $\bSigma_{\beta*}$.
In respect to Theorem \ref{theo31}, Condition A3 of \citet{CFLZ05},
ensures that $\max_{j} \sigma_j^2$ is bounded, therefore satisfying
$\max_j\sigma_j^2=O(u_n)$ for any positive sequence $u_n$ bounded away
from zero.

The following lemma controls the difference between the empirical
information matrix with
\[
\mathcal{I}_{\beta_1}=\int_0^T \bV(\bbeta_1,t) S_n^{(0)} \lambda
_0(t) \,dt
\]
and its population counterpart $\Sigma_{\beta_1}$, and plays a crucial
part in the theoretical developments of this section.
%
\begin{lemma} \label{le92a}
$\!\!\!$Assume that Conditions \ref{cond3} and \ref{cond6} hold. Then
$\sup_{\beta_1\in{\mathcal B}}\|{\mathcal I}_{\beta_1}\|_2=O_p(1)$,
$\|{\mathcal I}_{\beta_1^*}^{-1}\|_2=O_p(1)$
and
$\sup_{\beta_1\in{\mathcal B}}\|{\mathcal I}_{\beta_1}-\bSigma
_{\beta
_1}\|_2=o_p(1)$.
\end{lemma}
\begin{pf}
We prove the statement in the following three steps:

\begin{longlist}
\item
For any $s\times1$ vector function $\bfa (t)$ on $[0,\tau]$,
we have
\[
\biggl\|\int_0^{\tau}\bfa (t)\lambda_0(t) \,dt\biggr\|_2^2
\le\Lambda_0(\tau)\int_0^{\tau}\|\bfa (t)\|_2^2\lambda
_0(t) \,dt.
\]
In fact, by definition,
$\|\int_0^{\tau}\bfa (t)\lambda_0(t) \,dt\|_2^2
=\sum_{i=1}^s (\int_0^{\tau}a_i(t)\lambda_0(t) \,dt)^2$,
where~$a_i(t)$ is the $i$th component\vadjust{\goodbreak} function of $\bfa (t)$.
Using the H\"older inequality, we obtain that
\begin{eqnarray*}
\biggl\|\int_0^{\tau}\bfa (t)\lambda_0(t)
\,dt\biggr\|_2^2
&\le& \sum_{i=1}^s\Lambda_0(\tau) \int_0^{\tau}a_i^2(t)\lambda
_0(t)
\,dt\\
&=&\Lambda_0(\tau)\int_0^{\tau}\|\bfa (t)\|_2^2\lambda
_0(t) \,dt.
\end{eqnarray*}

\item For any matrix function $\bA(t)$ on $[0,\tau]$, we have
\[
\biggl\|\int_0^{\tau}\bA(t)\lambda_0(t) \,dt\biggr\|_2^2
\le\Lambda_0(\tau)\int_0^{\tau}\|\bA(t)\|_2^2\lambda_0(t) \,dt.
\]
In fact,
\begin{eqnarray*}
\biggl\|\int_0^{\tau}\bA(t)\lambda_0(t) \,dt\biggr\|_2^2
&=&\sup_{\|\bu\|_2=1} \biggl\|\biggl(\int_0^{\tau}\bA(t)\lambda
_0(t)
\,dt\biggr)\bu\biggr\|_2^2\\
&=&\sup_{\|\bu\|_2=1} \biggl\|\int_0^{\tau}\bfa _u(t)\lambda_0(t) \,dt
\biggr\|_2^2,
\end{eqnarray*}
where
$\bfa _u(t)=\bA(t)\bu$. Then
\begin{eqnarray*}
\int_0^{\tau}\|\bA(t)\|_2^2\lambda_0(t) \,dt
&=&\int_0^{\tau}\sup_{\|\bu\|=1}
\bu^T\bA(t)^{\otimes2}\bu\lambda_0(t) \,dt \\
&=&\int_0^{\tau}\sup_{\|\bu\|=1} \|\bfa _u(t)\|_2^2\lambda
_0(t) \,dt\\
&\ge&\sup_{\|\bu\|=1}\int_0^{\tau} \|\bfa _u(t)\|_2^2\lambda
_0(t) \,dt.
\end{eqnarray*}
Therefore, by (i), the result holds.

\item By definition, we have
\begin{eqnarray*}
{\mathcal I}_{\beta_1}-\bSigma_{\beta_1}
&=&\int_0^{\tau}\{\bV(\bbeta_1,t)-\bfv (\bbeta_1,t)\}
s^{(0)}(\bbeta
_1^*,t)\lambda_0(t) \,dt\\
&&{}+\int_0^{\tau}\bV(\bbeta_1,t)\bigl\{S_n^{(0)}(\bbeta
_1^*,t)-s^{(0)}(\bbeta
_1^*,t)\bigr\}\lambda_0(t) \,dt\\
&\equiv& \bA_{n1}(\bbeta_1)+\bA_{n2}(\bbeta_1).
\end{eqnarray*}
Using (ii), we obtain that
\[
\|\bA_{n1}(\bbeta_1)\|_2^2
\le \Lambda_0(\tau)\int_0^{\tau} \|\bV(\bbeta_1,t)-\bfv (\bbeta_1,t)\|
_2^2\bigl(s^{(0)}(\bbeta_1^*,t)\bigr)^2\lambda_0(t) \,dt.
\]
Then, by Condition \ref{cond3},
${\sup_{\beta_1\in{\mathcal B}}}\|\bA_{n1}(\bbeta_1)\|_2^2=o_p(1)$.
Similarly,
\[
\sup_{\beta_1\in{\mathcal B}}\|\bA_{n2}(\bbeta_1)\|_2^2=o_p(1).
\]
Therefore,
%
\begin{equation} \label{eqjb1}\qquad
\sup_{\beta_1\in{\mathcal B}}\|{\mathcal I}_{\beta_1}-\bSigma
_{\beta
_1}\|_2
\le\sup_{\beta_1\in{\mathcal B}}\|\bA_{n1}(\bbeta_1)\|_2
+\sup_{\beta_1\in{\mathcal B}}\|\bA_{n2}(\bbeta_1)\|_2=o_p(1).
\end{equation}
By Condition \ref{cond3}(ii), we have
\[
\sup_{\beta_1\in{\mathcal B}}\|\bSigma_{\beta_1}\|_2\le\int
_0^{\tau
}\sup_{\beta_1\in{\mathcal B},t\in[0,\tau]}\|\bfv (\bbeta
_1,u)\|
_2s^{(0)}(\bbeta_1^*,u) \,d\Lambda_0(u)=O_p(1).
\]
This combining with (\ref{eqjb1}) leads to
\[
\sup_{\beta_1\in{\mathcal B}}\|{\mathcal I}_{\beta_1}\|_2\le\sup
_{\beta
_1\in{\mathcal B}}\|\bSigma_{\beta_1}\|_2+\sup_{\beta_1\in
{\mathcal
B}}\|{\mathcal I}_{\beta_1}-\bSigma_{\beta_1}\|_2=O_p(1).
\]
Decompose ${\mathcal I}_{\beta_1^*}^{-1}$ as
\[
{\mathcal I}_{\beta_1^*}^{-1}=
\bSigma_{\beta_1^*}^{-1/2}\{I+\bSigma_{\beta_1^*}^{-1/2}({\mathcal
I}_{\beta_1^*}-\bSigma_{\beta_1^*})\bSigma_{\beta_1^*}^{-1/2}\}
^{-1}\bSigma_{\beta_1^*}^{-1/2}
\]
and let
${\mathcal A}=I+\bSigma_{\beta_1^*}^{-1/2}({\mathcal I}_{\beta
_1^*}-\bSigma_{\beta_1^*})\bSigma_{\beta_1^*}^{-1/2}$.
Then
${\mathcal I}_{\beta_1^*}^{-1}=\bSigma_{\beta_1^*}^{-1/2}{\mathcal
A}^{-1}\bSigma_{\beta_1^*}^{-1/2}$.\vspace*{1pt}
Using the Bauer--Fike inequality [\citet{B97}], we obtain that
\[
|\lambda({\mathcal A})-1|
\le\|\bSigma_{\beta_1^*}^{-1/2}({\mathcal I}_{\beta_1^*}-\bSigma
_{\beta
_1^*})\bSigma_{\beta_1^*}^{-1/2}\|_2
\le\|\bSigma_{\beta_1^*}^{-1/2}\|_2 \|{\mathcal I}_{\beta
_1^*}-\bSigma
_{\beta_1^*}\|_2
\|\bSigma_{\beta_1^*}^{-1/2}\|_2.
\]
Then by (\ref{eqjb1}) and Condition \ref{cond6},
$|\lambda({\mathcal A})-1|=o_p(1)$.
Hence,
$\lambda({\mathcal A}^{-1})=1+o_p(1)$.
Since ${\mathcal A}$ is symmetrical,
$\|{\mathcal A}^{-1}\|_2=O_p(1)$.
This together with Condition \ref{cond6}
yield that
$\|{\mathcal I}_{\beta_1^*}^{-1}\|_2\le
\|\bSigma_{\beta_1^*}^{-1/2}\|_2 \|{\mathcal A}^{-1}\|_2 \|\bSigma
_{\beta_1^*}^{-1/2}\|_2=O_p(1)$.
\end{longlist}
\upqed\end{pf}

The following tail condition is needed as a technicality in
establishing estimation loss results on the oracle estimator $\hat
\bbeta{}^{\mathbf{o}}$.
%
\begin{condition}\label{cond7}
$E\{\sup_{0\le t\le\tau}
Y(t)\|\bS(t)\|_{2}^2\exp(\bbeta_1^{*T}\bS(t))\}=O(s)$.
\end{condition}

For a fixed effective dimensionality $s$,
Condition \ref{cond7} is implied by the following condition from
\citet{AG82}:
%
\begin{equation}\label{condAG82}
E\Bigl\{\sup_{0\le t\le\tau,\beta_1\in{\mathcal B}}
Y(t)\|\bS(t)\|_{2}^2\exp(\bbeta_1^{T}\bS(t))\Bigr\}<\infty.
\end{equation}
However, we deal with diverging $s$, the above condition (\ref
{condAG82}) is obviously too tight to be satisfied.
For example, when all variables in $\mathcal{M}_*$ are bounded, we have
$\|\bS(t)\|_{2}^2=O(s)$. In general, if each $S_k(t)$ in $\bS(t)$
satisfies (\ref{condAG82}),
then Condition~\ref{cond7} holds. Now we are ready to state the result
on the existence of the biased oracle estimator.
%
\begin{theorem}[(Estimation loss)]\label{theo41}
Under Conditions \ref{cond1}, \ref{cond3} and \ref{cond6}, \ref{cond7},
with probability tending to one, there exists an oracle estimator
$\hat{\bbeta}{}^{\mathbf{o}}$ such that
\[
\|\hat{\bbeta}{}^{\mathbf{o}}-\bbeta^* \|_2 = O_P\bigl\{\sqrt{s}\bigl(n^{-1/2} +
\lambda_n\rho'(\beta_n^*)\bigr)\bigr\},
\]
where $\beta_n^*=\min\{|\beta_j^*|, j\in\mathcal{M}_* \}$ is the
minimum signal strength.\vadjust{\goodbreak}
\end{theorem}
\begin{pf}
Since $\hat{\bbeta}{}^{\mathbf{o}}_2=\bbeta_2^*=\mathbf{0}$, we only need
to consider the subvector in the first $s$ components, that is, we can
restrict our attention to the $s$-dimensional subspace
$\{ \bbeta_1\in\mathbb{R}^s\dvtx \bbeta_{\mathcal{M}_*^c}=0\}$.
It suffices to show that, for any $\varepsilon>0$, there exists a large
constant $B$ and
$\gamma_n=B\{\sqrt{s}(n^{-1/2} + \lambda_n \rho'(\beta_n^*)\}$
such that
\[
P\Bigl\{\sup_{\|\mathbf{u}\|_2=1} \mathcal{C}(\bbeta_1^*+\gamma_n \bu,
\mathbf{0}) < \mathcal{C}(\bbeta^*_1,\mathbf{0})\Bigr\}
\ge1-\varepsilon,
\]
when $n$ is big enough,
where for short
$\mathcal{C}({\bbeta})$ denotes
$\mathcal{C}({\bbeta},\tau)$,
and in particular
$\mathcal{C}({\bbeta}_1,\mathbf{0})$
represents $\mathcal{C}(({\bbeta}_1,\mathbf{0}),\tau)$.
This indicates that, with probability tending to one, there exists a
local maximizer such that
$\|\hat{\bbeta}{}^{{\mathbf o}}-\bbeta^*\|_2=O_p\{\sqrt{s}(n^{-1/2} +
\lambda_n \rho'(\beta_n^*))\}$.

Let $\bE_n^{(1)}(\bbeta_1,t)=\bE_n^{(1)}((\bbeta_1,\mathbf{0}),t)$,
$\mathcal{P}_n(\bbeta_1)=n\lambda_n\sum_{j=1}^s\rho(|\beta
_j|;\lambda_n)$
and
\[
U_n(\bbeta_1)=\partial{\mathcal L}(\bbeta_1)=\sum_{i=1}^n
\int_0^{\tau}\bigl\{\bS_i(t)-\bE_n^{(1)}(\bbeta_1,t)\bigr\} \,dN_i(t).
\]
By the Taylor expansion at $\gamma_n=0$,
%
\begin{eqnarray} \label{eq21}
&&\mathcal{C}(\bbeta_1^*+\gamma_n \bu, 0) -\mathcal{C}(\bbeta_1^*,
0)\nonumber\\
&&\qquad= \bu^TU_n(\bbeta_1^*) \gamma_n + 0.5\gamma_n^2 \bu^T \,\partial
U_n(\bbeta_1^*) \bu+ r_n({\bbeta}_1)\\
&&\qquad\quad{} -\mathcal{P}_n(\bbeta_1^*+\gamma_n \bu, 0)+\mathcal{P}_n(\bbeta
^*_1),\nonumber
\end{eqnarray}
where the remainder term $r_n({\bbeta}_1)$ is equal to
\[
\frac{1}{6}\sum_{j,k}
(\beta_{1j}-\beta^*_{1j})(\beta_{1k}-\beta^*_{1k})(\beta_{1\ell
}-\beta
^*_{1\ell})\,
\frac{\partial^2 U_{n\ell}(\bbeta_1)}{\partial\beta_{1j}\,\partial
\beta_{1k}}
\]
with $U_{n\ell}$ being the $\ell$th component of $U_{n}$
and $\bbeta_1$ lying between $\bbeta_1^*+\gamma_n\bu$ and~$\bbeta^*_1$.
By Lemma 2.2 in the supplementary material [\citet{BFJ11}] we have
$\|U_n(\bbeta_1^*)\|_2=O_p(\sqrt{ns})$.
It follows that
%
\begin{equation}\label{eqw0}
|\bu^T U_n(\bbeta_1^*) \gamma_n |=O_p\bigl(\sqrt{ns}\gamma_n\bigr).
\end{equation}
%
By simple\vspace*{1pt} decomposition, we have
$ \partial U_n(\bbeta_1^*) = -n({\mathcal I}_{\beta_1}+{\mathcal
W}_{\beta_1})$,
where $\mathcal{I}_{\beta_1}$ was defined in Lemma \ref{le93} and
${\mathcal W}_{\beta_1}= n^{-1} \int_0^\tau\bV(\bbeta_1,t) \,d \bar M(t)$.
Hence,
\begin{eqnarray*}
\gamma_n^2 \bu^T \,\partial U_n(\bbeta_1^*) \bu
&=&-n\gamma_n^2\{\bu^T(-n^{-1}\,\partial U_n(\bbeta_1^*))\bu\}\\
&=&-n\gamma_n^2\{\bu^T\bSigma_{\beta_1^*}\bu+\bu^T[({\mathcal
I}_{\beta
_1^*}-{\bSigma}_{\beta_1^*})+{\mathcal W}_{\beta_1^*}]\bu\}.
\end{eqnarray*}
By Lemma 2.3 in the supplementary material [\citet{BFJ11}] and Lemma
\ref{le92a},
\[
\|({\mathcal I}_{\beta_1^*}-{\bSigma}_{\beta_1^*})+{\mathcal
W}_{\beta
_1^*}\|_2
\le\|{\mathcal I}_{\beta_1^*}-{\bSigma}_{\beta_1^*}\|_2+\|{\mathcal
W}_{\beta_1^*}\|_2
=o_p(1).
\]
Therefore, by Condition \ref{cond6},
there exists a constant $c>0$ such that
%
\begin{equation}\label{eqw2}
\gamma_n^2 \bu^T \,\partial U_n(\bbeta_1^*) \bu\leq-cn\gamma_n^2\bigl(1+o_p(1)\bigr).
\end{equation}
Since
$\|\bbeta_1-\bbeta_1^*\|_2\le\gamma_n$
and the average of i.i.d. terms,
$n^{-1}\,\frac{\partial^2 U_{n\ell}(\bbeta_1)}{\partial\beta
_{1j}\,\partial
\beta_{1k}}$, is of order $O_p(1)$,
we have
$
r_n(\bbeta_1)=O_p(n\gamma_n^3).
$
By concavity\vspace*{1pt} of $\rho$ and decreasing property of $\rho'$ from
Condition \ref{cond1},
\begin{eqnarray*}
|\mathcal{P}_n(\bbeta_1^*+\gamma_n \bu, 0)-\mathcal{P}_n(\bbeta
^*_1) |
&=& n\lambda_n \sum_{j=1}^s \bigl|\rho(|\beta_j^*+\gamma_nu_j|;\lambda
_n)-\rho(|\beta_j^*|;\lambda_n)\bigr|\\
&\leq& n \lambda_n \gamma_n \|\brho'_0(\beta_n^*)\| \|\bu\|_2\bigl(1+o_p(1)\bigr),
\end{eqnarray*}
where
$\beta_n^*$ is the minimal signal length and
$\brho'_0(\cdot)$ is the subvector of $\brho(\cdot)$, consisting of its
first $s$ elements.
Then
%
\begin{equation}\label{eq23}
|\mathcal{P}_n(\bbeta_1^*+\gamma_n \bu, 0)-\mathcal{P}_n(\bbeta
^*_1) |
=O_p\bigl(n\lambda_n\sqrt{s}\gamma_n\rho'(\beta_n^*) \bigr).
\end{equation}
Combining (\ref{eq21})--(\ref{eq23}) leads to
\[
\mathcal{C}(\bbeta_1^*+\gamma_n \bu, 0) -\mathcal{C}(\bbeta_1^*, 0)
< n \gamma_n
\bigl\{O_p\bigl(\sqrt{s/n} + \sqrt{s} \lambda_n\rho'(\beta_n^*)\bigr) - c \gamma
_n\bigl(1+o_p(1)\bigr) \bigr\},
\]
where with probability tending to one, the RHS is smaller then zero
when $\gamma_n=B(\sqrt{s/n} + s \lambda_n \rho'(\beta_n^*))$ for a
sufficiently large $B$.
\end{pf}

A simple corollary of this theorem is that the $L_{1},L_\infty$
estimation losses of the oracle estimator are bounded by $ s(n^{-1/2} +
\lambda_n\rho'(\beta_n^*))\}$ and by $\sqrt{s}(n^{-1/2} + \lambda
_n\rho
'(\beta_n^*))$, respectively. Hence, $L_1$ loss can have a chance to be
close to zero only if the sparsity parameter $\alpha<1/2$, whereas
$L_{\infty}$ loss will converge to zero with no restrictions on
$\alpha$.

To make the bias in the penalized estimation negligible, $\rho'(\beta
_n^*)$ needs to converge to zero at a specific rate controlled by the
next condition.
%
\begin{condition}\label{cond5}
The regularization\vspace*{1pt} parameter $\lambda_n$ satisfies that
$\sqrt{s} \lambda_n \rho'(\beta_n^*;\allowbreak\lambda_n) \to0$ and
$\lambda_n\gg n^{-0.5+(0.5\alpha+\alpha_1-1)_++\alpha_2}$,
where $\alpha_1$ is defined in Condition~\ref{cond4}, and $\alpha_2$
is a positive constant.
\end{condition}

Condition \ref{cond5} regulates the behavior of the regularization
parameter $\lambda_n$ around $0$ and $\infty$. From the result of
Theorem \ref{theo41}, we see that for different penalties, the ``extra
term'' $\sqrt{s} \lambda_n \rho'(\beta_n^*;\lambda_n) $ in the $L_2$
estimation loss will require either extra conditions on the $\lambda_n$
or extra conditions on the minimum signal strength $\beta_n^*$ (see
Theorems \ref{theo42}--\ref{cor3} for further details) and can govern
estimation efficiency of the penalized estimators.
%
\begin{condition}\label{cond5a}
$\!\!\!$Let $\kappa_0=\max_{\delta\in{\mathcal N}_0}\kappa(\rho,\delta)$,
where
${\mathcal N}_0=\{\delta\in R^s\dvtx \|\delta-\bbeta_1^*\|_{\infty
}\le
\beta_n^*\}$.
Assume that
$\lambda_n$ and $\beta_n^*$ satisfy that
(i) $\beta_n^*\gg\sqrt{s}(n^{-1/2}+\lambda_n\rho'(\beta_n^*))$
and
(ii) $\lambda_{\min}(\bSigma_{\beta_1^*})>\lambda_n\kappa_0$.
\end{condition}

Condition \ref{cond5a}(i) is employed to make $\hat{\bbeta
}{}^{\mathbf o}_1$ fall in ${\mathcal N}_0$
with probability tending to one.
For LASSO, since $\rho'(\beta_n^*)=1$, it means\vspace*{1pt} that
$\beta_n^*\gg\sqrt{s}\lambda_n$. By Condition~\ref{cond5},
it reduces to $\beta_n^*\gg\sqrt{s}n^{-0.5+(0.5\alpha+\alpha
_1-1)_+\alpha_2}$.
For SCAD,
if $\beta_n^*\gg\lambda_n$, then $\rho'(\beta_n^*)=0$ when $n$ is
large enough
and hence it requires that
$\beta_n^*\gg\sqrt{s}n^{-0.5}$.
Therefore, Condition \ref{cond5a}(i) is less restrictive for
SCAD-like penalties.
Condition \ref{cond5a}(ii) is used to ensure the condition in (\ref
{eq7}) holds with probability tending to one (see the proof of Theorem
\ref{theo42}). It always holds when $\kappa_0 = 0$ (e.g., for the
LASSO penalty) and is satisfied for the SCAD type of penalty when
$\beta
_n^*\gg\lambda_n$.
%
\begin{condition}\label{cond4} For $\alpha_1 >0 $
and $0<C<\infty$,
\[
{\sup_{0\leq t\leq\tau} \sup_{\mathbf{v}_1 \in\mathcal
{B}(\bolds
\beta_1^*, \beta_n^*)}}\|
\tilde{\bV}(t, \bfv )
\|_{2,\infty} = \min\biggl( C\frac{\rho'(0+)}{\rho'(\beta_n^*)},
O_p( n^{\alpha_1}) \biggr),
\]
where\vspace*{1pt}
$\mathcal{B}(\bbeta_1^*,\beta_n^*)$ is an $s$-dimensional ball
centered at $\bbeta_1^*$ with radius $\beta_n^*$,
for $\bfv =(\bfv _1^T,\mathbf{0}^T)^T$,
\[
\tilde{\bV}(t,\bfv )=\frac{S^{(0)}_{n}(\bfv ,t)
S^{(2)}_{n21}(\bfv ,t)
- {{S^{(1)}_{n2}}(\bfv ,t)} {(S^{(1)}_{n1}}(\bfv ,t))^T
}{\{S^{(0)}_{n}(\bfv
,t)\}^{2} } \in\mathbb{R}^{(p-s)\times s}
\]
and
$\|\tilde{\bV}(t,\bfv )\|_{2,\infty}={\max_{\|\bfx \|
_2=1}}\|\tilde{\bV}(t,\bfv
)\bfx \|_{\infty}$.
\end{condition}

As noted in Fleming and Harrington [(\citeyear{FH05}), page 149],
$\tilde{\bV} (\bbeta,t)$ is an empirical covariance matrix of $\bX
_i(t)$ computed with
weights proportional to $Y_i(t)\times\allowbreak\exp\{\bbeta^{T}\bX_i(t)\}$. Hence,
$\tilde{\bV}$ is the empirical covariance matrix between the important
variables $\bS_i(t)$ and
unimportant variables $\bQ_i(t)$.
Condition \ref{cond4} controls the uniform growth rate of the norm of
these covariance matrices, a
notion of weak correlation between $\bS_i(t)$ and $\bQ_i(t)$.
For the $L_1$ penalty, $\rho'(\beta_n^*)=1$, and Condition \ref
{cond4} becomes
a version of ``strong irrepresentable'' condition [\citet{ZY06}] for
censored data. It is
very stringent as the right-hand side has to be bounded by $O(1)$. On
the other hand for the SCAD penalty, if $\beta_n^*\gg\lambda_n$, then
$\rho'(\beta_n^*)=0$ when $n$ is large enough. Therefore,
Condition \ref
{cond4} is significantly relaxed to $O(n^{\alpha_1})$.
In general, when a folded concave penalty is employed, the upper bound
on the right-hand side in Condition \ref{cond4} can grow to infinity
at polynomial rate. This was also noted in the work of \citet{FL10} in
the context of generalized linear models.
%
\begin{theorem}[(Strong oracle)] \label{theo42}
Let the oracle estimator $\hat{\bolds\beta}{}^{\mathbf{o}}$ be a
local~ma\-ximizer of $\mathcal{C}(\bbeta_1,\tau)$ given by Theorem
\ref{theo41}.
If
$\max_j(\sigma_j^2)=O(n^{(0.5\alpha+\alpha_1-1)_++\alpha_2})$, and
Conditions \ref{cond1}--\ref{cond4} hold, then with
probability tending to one,
there exists a local maximizer $\hat\bbeta$ of $\mathcal{C}(\bbeta
,\tau
)$ such that
\[
P(\hat{\bolds\beta} = \hat{\bolds\beta}{}^{\mathbf
{o}})\geq
1-c_0(p-s) \exp\bigl\{ - c_1{n^{(0.5\alpha+\alpha_1-1)_++\alpha_2}}
\bigr\},
\]
where
$c_0$ and $c_1$ are positive constants.\vadjust{\goodbreak}
\end{theorem}
\begin{pf}
It suffices to show that $\hat{\bbeta}{}^{\mathbf{o}}$ is a local
maximizer of $\mathcal{C}(\bbeta,\tau)$ on a set
$\Omega_n$ which has a probability tending to one.
By Theorem \ref{theo21}, we need to show that, with probability
tending to one,
$\hat{\bbeta}{}^{\mathbf{o}}$ satisfies (\ref{eq5})--(\ref{eq7}).
Since $\hat{\bbeta}{}^{\mathbf{o}}$ already satisfies (\ref{eq5}) by
definition, we are left to check (\ref{eq6})
and (\ref{eq7}).

Define
$\Omega_n=\{\bxi\dvtx \|\bxi_{\mathcal{M}_*^c} \|_{\infty} \leq
\sqrt
{n}u_n \}$
for some diverging sequence $u_n$ to be chosen later,
where $\bxi_{\mathcal{M}_*^c}$ is the subvector of $\bxi$ with indices
in $\mathcal{M}_*^c$.
By Theorem \ref{theo31}, there exist positive constants $c_0$ and
$c_1$ such that
\[
P \bigl( |\xi_j|>\sqrt{n}u_n \bigr)
\leq c_0 \exp\{-c_1u_n\}
\]
uniformly over $j$.
Then using the Bonferroni union bound, we obtain that
%
\begin{eqnarray}\label{temp2}
P(\Omega_n) &\geq& 1-\sum_{j \in\mathcal{M}_*^c}P \bigl( |\xi
_j|> \sqrt
{n}u_n \bigr)\nonumber\\[-8pt]\\[-8pt]
&\geq& 1- c_0(p-s)e^{-c_1 u_n} \to1\qquad \mbox{as } n\to
\infty,\nonumber
\end{eqnarray}
where $u_n$ can be chosen later to make $(p-s)e^{-c_1u_n}\to0$.
We now check if~(\ref{eq6}) holds for $\hat{\bbeta}{}^{\mathbf{o}}$ on
the set $\Omega_n$.
Denote by $\brho'_{\mathcal{M}_*^c}$ the subvector of $\brho'(|\hat
{\bbeta}{}^{\mathbf{o}}|)$ with indexes in ${\mathcal{M}_*^c}$.
Let
$\gamma(\bbeta)=\int_0^t S^{(1)}_{n}(\bbeta,u)/S^{(0)}_{n}(\bbeta
,u)\,
d\bar{N}(u)$
and
\[
\mathbf{z}(\hat{\beta}{}^{\mathbf{o}})
= \sum_{i=1}^n
\int_0^{\tau}\bigl\{
\bQ_i(t)-\bE_n^{(2)}(\hat{\bbeta}{}^{\mathbf{o}},t)\bigr\} \,dN_i(t),
\]
where
$\bE_n^{(2)}(\bbeta,t)=S_{n2}^{(1)}(\bbeta,t)/S_n^{(0)}(\bbeta,t)$.
Then by Condition \ref{cond1}, we have
%
\begin{eqnarray} \label{eq28}
\|\mathbf{z}(\hat{\bbeta}{}^{\mathbf{o}})\|_{\infty}
& \leq& \| \bxi_{\mathcal{M}_*^c} \|_{\infty}
+ \| \gamma_{\mathcal{M}_*^c}(\bbeta^*)
-\gamma_{\mathcal{M}_*^c}(\hat{\bbeta}{}^{\mathbf{o}})
\|_\infty\nonumber\\
& = & O\biggl( \sqrt{n}u_n
+ \biggl\| \int_0^t \tilde{V}(u, \bfv _1) (\hat\bbeta
{}^{\mathbf{o}}_{1} -
\bbeta_1^*)\,d\bar N(u)
\biggr\|_{\infty} \biggr)
\\
&=& O\Bigl(\sqrt{n} u_n
+ {\sup_{0\leq u\leq\tau} \sup_{\mathbf{v}_1\in\mathcal
{B}(\bolds
\beta_1,\beta_n^*)}} \| \tilde{V}(u, \bfv _1)
\|_{2,\infty} \| \hat\bbeta{}^{\mathbf{o}}_{1} - \bbeta_1^*
\|_2 \Bigr),\nonumber
\end{eqnarray}
where $\bfv =(\bfv _1^T,\mathbf{0}^T)^T$ with $
\bfv_1$ being between $\bbeta
_1^*$ and $\bbeta{}^{\mathbf{o}}_{1}$, and $\tilde{V}(u,{\bfv }_1)$ is
defined in Condition \ref{cond4}.
By Theorem \ref{theo41} and Condition \ref{cond4},
we obtain that
$(n \lambda_n\times\rho'(0+))^{-1} \|\mathbf{z}(\hat{\bbeta}{}^{\mathbf
{o}})\|
_{\infty}$ is bounded by
\begin{eqnarray*}
&& n^{-1}\lambda_n^{-1}O_p\Bigl\{\sqrt{n} u_n
+ {\sup_{0\leq u\leq1}}\| \tilde{V}(u, \bfv _1)
\|_{2,\infty}\sqrt{s}\bigl(n^{-1/2}+\lambda_n\rho'(\beta_n^*)
\bigr) \Bigr\}
\\
&&\qquad= O_p\{n^{-1/2}\lambda_n^{-1} (u_n + n^{0.5\alpha+\alpha_1 -1})
+ n^{-1+0.5\alpha}\rho'(0+) \} \to0,
\end{eqnarray*}
if we take
$u_n=n^{(0.5\alpha+\alpha_1-1)_{+}+\alpha_2}$
and
$\lambda_n\gg n^{-0.5+(0.5\alpha+\alpha_1-1)_{+}+\alpha_2}$.
Therefore,~(\ref{eq6})~holds on $\Omega_n$.
Once $\delta<(0.5\alpha+\alpha_1-1)_{+}+\alpha_2$, (\ref{temp2}) holds.

We are now left to show that (\ref{eq7}) holds with probability tending
to one, that is,
$\lambda_{\min}\{n^{-1}\int_0^{\tau}\bV(\hat{\bbeta
}{}^{\mathbf
{o}},t) \,d\bar{N}(t)\}
> \lambda_n\kappa(\rho,\hat{\bbeta}{}^{{\mathbf o}}_1)$,
which is guaranteed by Condition~\ref{cond5a}. In fact, by Theorem
\ref{theo41} and Condition \ref{cond5a}(i), with probability tending to
one, $\hat{\bbeta}{}^{{\mathbf o}}_1$~falls in ${\mathcal N}_0$ as
$n\to\infty$, so that $\kappa(\rho,\hat{\bbeta}{}^{\mathbf
o}_1)\le\kappa_0$. Hence, by Condition \ref{cond5a}(ii), with
probability tending to one,
%
\begin{equation}\label{addeqs1}
\lambda_{\min}(\bSigma_{\beta_1^*})> \lambda_n\kappa(\rho,\hat
{\bbeta
}{}^{\mathbf o}_1).
\end{equation}
Recall that
\begin{eqnarray*}
n^{-1}\int_0^{\tau}\bV(\hat{\bbeta}{}^{\mathbf o}_1,t) \,d\bar{N}(t)
&=&{\mathcal I}_{\hat{\beta}{}^{\mathbf o}_1}+{\mathcal W}_{\hat{\beta
}{}^{\mathbf o}_1}
=\bSigma_{\hat{\beta}{}^{\mathbf o}_1}
+({\mathcal I}_{\hat{\beta}{}^{\mathbf o}_1}-\bSigma_{\hat{\beta
}{}^{\mathbf o}_1})
+{\mathcal W}_{\hat{\beta}{}^{\mathbf o}_1}.
\end{eqnarray*}
By Theorem \ref{theo41}, as $n\to\infty$,
$\hat{\bbeta}{}^{\mathbf o}_1\in{\mathcal B}$ with probability tending
to one.
This combining with
Lemma \ref{le92a} and Lemma 2.3 in the supplementary material [\citet{BFJ11}] leads to
\[
n^{-1}\int_0^{\tau}\bV(\hat{\bbeta}{}^{\mathbf o}_1,t) \,d\bar{N}(t)
=\bSigma_{\hat{\beta}{}^{\mathbf o}_1}+E,
\]
where $\|E\|_2=o_p(1)$.
By Condition \ref{cond3}(i), (iii),
with probability tending to one,
\[
\|\bSigma_{\hat{\beta}{}^{\mathbf o}_1}
-\bSigma_{\beta_1^*}\|_2=o_p(1).
\]
Let\vspace*{1pt} $E^*=n^{-1}\int_0^{\tau}\bV(\hat{\bbeta}{}^{\mathbf o}_1,t)
\,d\bar{N}(t)
-\bSigma_{\beta_1^*}$. Then
$\|E^*\|_2=o_p(1)$.
Using Weyl's pertubation theorem [\citet{B97}],
we obtain that
\[
\min_{1\le k\le s}\biggl|\lambda_k \biggl\{n^{-1}\int_0^{\tau}\bV(\hat
{\bbeta
}{}^{\mathbf o}_1,t) \,d\bar{N}(t)\biggr\}
-\lambda_k(\bSigma_{\beta_1^*})\biggr|\le\|E^*\|_2,
\]
where
$\lambda_k(\bSigma_{\beta_1^*})$ is the $k$th largest eigenvalue of
$\bSigma_{\beta_1^*}$.
Therefore,
\[
\lambda_{\min}\biggl\{n^{-1}\int_0^{\tau}\bV(\hat{\bbeta
}{}^{\mathbf
o}_1,t) \,d\bar{N}(t)\biggr\}
=\lambda_{\min}(\bSigma_{\beta_1^*})+o_p(1).
\]
This combining with (\ref{addeqs1}) yields that with probability
tending to one
\[
\lambda_{\min}\biggl\{n^{-1}\int_0^{\tau}\bV(\hat{\bbeta
}{}^{\mathbf
o}_1,t) \,d\bar{N}(t)\biggr\}
>\lambda_n\kappa(\rho,\hat{\bbeta}{}^{\mathbf o}_1).
\]
\upqed
\end{pf}

The theorem becomes nontrivial if $\delta<(0.5\alpha+\alpha
_1-1)_++\alpha_2$, since $\log p = O(n^\delta)$.
Apart from the work of \citet{BFW09}, no formal work explicitly relates
the oracle property and the full and effective dimensionalities.
Theorem \ref{theo42} shows that $\hat{\bolds\beta}$ becomes
the biased
oracle with probability tending to one exponentially fast.
Then combining Theorems \ref{theo41} and \ref{theo42} leads to the
following $L_2$ estimation loss:
%
\begin{equation}\label{eqjj1}
\|\hat{\bbeta}_1-\bbeta^*_1 \|_2 = O_P\bigl\{\sqrt{s}\bigl(n^{-1/2} +\lambda
_n\rho
'(\beta_n^*)\bigr)\bigr\}.
\end{equation}
%
This theorem tells us that the resulting estimator behaves as if the
true set of ``important variables''\vadjust{\goodbreak} (i.e., as oracle estimator) were
known with probability converging to 1 as both $p$ and
$n$ go to $\infty$. The previous notions of oracle were that the
estimator behaves like the oracle
rather than an actual oracle itself. Classical oracle property of
\citet
{FL02} or sign consistency of \citet{BRT09} are both corollaries of
this result. In this sense Theorem \ref{theo42} introduces a tighter
notion of an oracle property.
It was first mentioned in \citet{KCO08} for the SCAD estimator of the
linear model with polynomial dimensionality
and then extended by \citet{BFW09} to the penalized M-estimators under
the ultra-high dimensionality setting. Extending their work to Cox's
model was exceptionally challenging because of martingale and censoring
structures.
%
\begin{theorem}[(LASSO)]\label{cor2}
Under Conditions \ref{cond3}--\ref{cond7}, if
$\max_j(\sigma_j^2)=O(n^{\alpha_2})$, $\sqrt{s}\lambda_n \to0$,
$\lambda_n\gg n^{-0.5+\alpha_2}$ and
\[
{\sup_{0\leq t\leq\tau} \sup_{\mathbf{v}_1 \in\mathcal
{B}(\bolds
\beta_1^*, \beta_n^*)}}\|
\tilde{\bV}(t, \bfv )
\|_{2,\infty} = O_p(1),
\]
then the result in Theorem \ref{theo42} holds for \textup{LASSO}
estimator with probability being at least $1- c_0(p-s)\exp\{ -
c_1n^{\alpha_2}\}$. Furthermore,
\[
\|\hat\bbeta_1 -\bbeta_1^*\|_2=O_P\bigl(\sqrt{s}\lambda_n\bigr).
\]
\end{theorem}

The proof of this theorem is relegated to the supplementary material [\citet{BFJ11}].
For the LASSO, the rate of convergence for nonvanishing components is
dominated by the bias term $\lambda_n \gg n^{-1/2}$. In addition, since
$s=n^{\alpha}$,
the condition $\sqrt{s}\lambda_n \to0$
indicates that $\alpha<1-2\alpha_2$, where $\alpha_2 \in[0,1/2)$. That
is, the bigger is $\alpha_2$, and the smaller sparsity dimension $s$
can be recovered using LASSO. Moreover, LASSO with $\alpha_2<1/2$
requires $p\ll\exp\{c_1n^{\alpha_2}\}$ to achieve the strong oracle property.
Hence, as $p$ (or $\alpha_2$) gets bigger, $s$ (or $\alpha$) should
get smaller.
This means that, as data dimensionality gets higher, recoverable
problems get sparser.
This is a new discovery and has not been documented in the literature.
On the other hand for folded concave penalties, faster rates of
convergence are obtained with fewer restrictions on $p$ and $s$.
This can be seen from the following result, which is a straightforward
corollary of Theorem \ref{theo42} and whose proof is left for the
supplementary material [\citet{BFJ11}].
%
\begin{theorem}[(SCAD)] \label{cor3}
Under Conditions \ref{cond1}--\ref{cond7}, if
$\beta_n^* \gg\lambda_n$,
$\max_j(\sigma_j^2)=O(n^{(0.5\alpha+\alpha_1-1)_++\alpha_2})$,
$\lambda_n\gg n^{-0.5+(0.5\alpha+\alpha_1-1)_++\alpha_2}$ and
\[
{\sup_{0\leq t\leq\tau} \sup_{\mathbf{v}_1 \in\mathcal
{B}(\bolds
\beta_1^*, \beta_n^*)}}\|
\tilde{\bV}(t, \bfv )
\|_{2,\infty} = O_p(n^{\alpha_1}),
\]
then the\vspace*{1pt} result in Theorem \ref{theo42} holds for \textup{SCAD} estimator
with probability being at least $1- c_0(p-s)\exp\{ - c_1n^{(0.5\alpha
+\alpha_1-1)_++\alpha_2 }\}$. Furthermore,
\[
\|\hat\bbeta_1 -\bbeta_1^*\|_2=O_P\bigl(\sqrt{s/n}\bigr).
\]
\end{theorem}


Note that the proof of Theorem \ref{theo42} shows $\hat{\bbeta
}_2=\mathbf{0}$ on a set whose probability measure is going to one
exponentially fast. For statistical inference about~$\bbeta$,
asymptotic properties of $\hat{\bbeta}_1$ are needed to be explored. To
be able to construct confidence intervals of ${\bbeta}_1$ we need to
derive its asymptotic distribution. This was done in \citet{FL02} for
fixed $p$
and in \citet{CFLZ05} for $p=o(n^{1/4})$.
Here, we allow $p$ to diverge at exponential rate~$O(\exp\{n^{\delta
}\}
)$ and the effective dimensionality $s$ to diverge at rate of $o(n^{1/3})$.
To the best of our knowledge there is no work available for such a setting.
Extending the previous work to such a NP-dimensional setting is not
trivial and requires complicated eigenvalue results. Moreover, the
large deviation result in Section \ref{sec3}, the strong oracle result
in Theorem \ref{theo42} and Lemmas~2.1--2.3 in the
supplementary material [\citet{BFJ11}] are essential for establishing the desired
asymptotics. Moreover, the following Lemma \ref{le93} is an important
extension of the classical asymptotic Taylor expansion results when the
number of parameters is diverging with the sample size.
%
\begin{lemma} \label{le93}
For any $s\times1$ unit vector $\bfb _n$, let
\[
\phi_n=\bfb _n^T\bSigma_{\beta_1^*}^{1/2}(-n^{-1}\,\partial
U_n(\bbeta
_1^*))^{-1}n^{-1/2}U_n(\bbeta_1^*)
\]
and
\[
\phi_{n1}=\bfb _n^T\bSigma_{\beta_1^*}^{-1/2}
n^{-1/2}U_n(\bbeta_1^*).
\]
If Conditions \ref{cond3}, \ref{cond6} and \ref{cond7} hold and if
$s=o(n^{1/3})$, then
$\phi_n=\phi_{n1}+o_p(1)$.
\end{lemma}
\begin{pf}
Let\vspace*{1pt} ${\mathcal B}=I+{\mathcal I}_{\beta
_1^*}^{-1/2}{\mathcal
W}_{\beta_1^*}{\mathcal I}_{\beta_1^*}^{-1/2}$,
where $I$ is an $s\times s$ identity matrix.
Using the Bauer--Fike inequality [\citet{B97}], we obtain that
\[
|\lambda({\mathcal B})-1|\le\|{\mathcal I}_{\beta
_1^*}^{-1/2}{\mathcal W}_{\beta_1^*}{\mathcal I}_{\beta
_1^*}^{-1/2}\|_2.
\]
Then by the H\"older inequality we have
$|\lambda({\mathcal B})-1|\le\|{\mathcal I}_{\beta
_1^*}^{-1/2}\|_2^2\|{\mathcal W}_{\beta_1^*}\|_2$.
Applying Condition \ref{cond6} and Lemma \ref{le92a} and Lemma
2.3 of the supplementary material [\citet{BFJ11}], we establish that
%
\begin{equation}\label{eqjjca2}
\lambda({\mathcal B})=1+O_p\bigl(s/\sqrt{n}\bigr)
\end{equation}
uniformly for all eigenvalues of ${\mathcal B}$.
Note that
\[
(-n^{-1}\,\partial U_n(\bbeta_1^*))^{-1}
=({\mathcal I}_{\beta_1^*}+{\mathcal W}_{\beta_1^*})^{-1}
={\mathcal I}_{\beta_1^*}^{-1}
-{\mathcal I}_{\beta_1^*}^{-1/2}
\{I-{\mathcal B}^{-1}\}{\mathcal I}_{\beta_1^*}^{-1/2}.
\]
It follows that
\begin{eqnarray*}
\phi_n
&=&\bfb _n^T\bSigma_{\beta_1^*}^{1/2}{\mathcal
I}_{\beta
_1^*}^{-1} n^{-1/2}U_n(\bbeta_1^*)
-\bfb _n^T\bSigma_{\beta_1^*}^{1/2}{\mathcal I}_{
\beta
_1^*}^{-1/2}\{I-{\mathcal B}^{-1}\}{\mathcal I}_{\beta
_1^*}^{-1/2}n^{-1/2}U_n(\bbeta_1^*)\\
&\equiv& \phi_{n1}-\phi_{n2}.
\end{eqnarray*}
Since $I-{\mathcal B}^{-1}$ is symmetrical,
$r_{\sigma}(I-{\mathcal B}^{-1})=\|I-{\mathcal B}^{-1}\|_2$.
Recall that $\|\bfb _n\|_2=1$;
it follows that
\[
|\phi_{n2}|
\le r_{\sigma}(I-{\mathcal B}^{-1}) \|\bSigma_{\beta_1^*}^{1/2}\|_2
\|{\mathcal I}_{\beta_1^*}^{-1/2}\|_2^2
\|n^{-1/2}U_n(\bbeta_1^*)\|_2.
\]
By\vspace*{1pt} Condition \ref{cond6}, $\|\bSigma_{\beta_1^*}^{1/2}\|_2=O_p(1)$.
From Lemma \ref{le92a}, we have
$\|{\mathcal I}_{\beta_1^*}^{-1/2}\|_2=O_p(1)$.
By Lemma 2.2 in the supplementary material [\citet{BFJ11}],
$\|n^{-1/2}U_n(\bbeta_1^*)\|_2=O_p(\sqrt{s})$.
Therefore,
%
\begin{equation}\label{eqjjca3}
|\phi_{n2}|=r_{\sigma}(I-{\mathcal B}^{-1})O_p\bigl(\sqrt{s}\bigr).
\end{equation}
By definition, it is easy to see that
\[
r_{\sigma}(I-{\mathcal B}^{-1})=\max\{|1-\lambda|\dvtx \lambda\in
\sigma
({\mathcal B}^{-1})\}
=\max\{|1-\lambda^{-1}|\dvtx \lambda\in\sigma({\mathcal B})\},
\]
which, combined with (\ref{eqjjca2}), leads to
$r_{\sigma}(I-{\mathcal B}^{-1})=O_p(s/\sqrt{n})$.
This together with (\ref{eqjjca3}) yields that $\phi_{n2}=O_p(\sqrt
{s^3/n})=o_p(1)$, if $s=o(n^{1/3})$.
Hence, $\phi_n=\phi_{n1}+o_p(1)$.
\end{pf}

With the Lemma \ref{le93} and technical lemmas presented in the
supplementary material [\citet{BFJ11}] we are ready to state the results on the
asymptotic behavior of the penalized estimator. Detailed proof is
included in the supplementary material [\citet{BFJ11}].
%
\begin{theorem}\label{theo61}
Under Conditions \ref{cond1}--\ref{cond4}, and for $\lambda_n\rho
'(\beta_n^*)=o( (sn)^{-1/2})$ for any $s\times1$ unit vector $
\bfb_n$,
if $s=o(n^{1/3})$, the penalized partial likelihood estimator $\hat
{\bbeta}_1$ from (\ref{eqjj1}) satisfies
\[
\sqrt{n}\mathbf{b}_n^T{\bSigma}^{1/2}_{\beta_1^*} (\hat{\bbeta
}_1-\bbeta
_1^*) \to\mathcal{N}(0, 1).
\]
\end{theorem}

Theorems \ref{theo42} and \ref{theo61} claim that $\hat{\bbeta}$
enjoys model selection consistency
and achieves the information bound mimicking that of the oracle estimator
$\hat{\bbeta}{}^{\bo}$.

\section{Iterative coordinate ascent algorithm (ICA)}\label{sec7}

Coordinate-wise algorithms are especially attractive for $p\gg n$ and
have been previously introduced for penalized least-squares with the
$L_q$-penalty by \citet{DDD04}, \citet{FHHT07}, \citet{WL08} and for
generalized linear models with the folded concave penalty by \citet
{FHT10} and \citet{FL10}.
By Condition \ref{cond1} and Proposition 2.7.1 in \citet{B03}, the
coordinate-wise
maximization algorithm in each iteration provides limits that are
stationary points
of the overall optimization (\ref{eq4}). Therefore, each output of ICA
algorithm will give a stationary point. We will
adapt the algorithm in \citet{FL10} to the censored data.

First, let us, with slight abuse in notation, denote by $Q_n(\bbeta
)=L_n(\bbeta)-P_n(\bbeta)$, where $L_n(\cdot)$ and $P_n(\cdot)$ stand
for the loss and penalty parts, respectively. Let $l_{n}(\bbeta,\zeta
,j)$ be the \textit{partial quadratic approximation} of $L_n(\bbeta)$ at
$\zeta\in R^p$ along the $j$th coordinate, where $\{\beta_{k} = \zeta
_k, k\neq j\}$ are held fixed, but~$\beta_j$ is allowed to vary
\[
q_{n}(\bbeta_j,\zeta,j)=l_{n}(\bbeta,\zeta,j)- n p_{\lambda
_n}(|\beta_j|).
\]
Because of the complex likelihood function we need an additional loop
to compute the partial quadratic approximation.

This penalized quadratic optimization problem can be solved
analytically, avoiding the challenges of nonconcave optimization.
It updates each coordinate if the maximizer of the penalized
univariate optimization strictly increases the objective function
$Q_n(\bbeta)$ and if it satisfies $\{j\dvtx|z_j|>\rho'(0+)\}$. The
algorithm stops when two values of the objective function $Q_n(\bbeta)$
are not different by more than~$10^{-8}$, say. Details of the algorithm
are presented in the supplementary material [\citet{BFJ11}].


\subsection{Simulated examples} \label{subsec72}

To show good model selection and estimation properties of the proposed
methodology,
we simulated 100 standard Toeplitz ensembles of size 100 with
population correlation $\rho(X_i, X_j)=\rho^{|i-j|}$ with $\rho$
ranging from $0.25$, $0.5$, $0.75$ and $0.9$. The distribution of
censoring time $C$ is exponential with mean $U * \exp\{\bX_i^T \bbeta
\}
$, where $U$ is randomly generated from uniform distribution over
$[1,3]$ for each simulated data set. This censoring was used in \citet
{FL02}, which makes about 30\% of the data censored. The full and
effective dimensionalities of the true parameter $\bbeta$ are taken as
$\{100,4\}$, $\{1\mbox{,}000,4\}$, $\{5\mbox{,}000,4\}$ and $\{1\mbox{,}000, 25\}$,
respectively, with values $\pm1$ randomly placed (the rest is set as
zero). The penalties employed are LASSO [\citet{T96}], SCAD [\citet{FL01}],
SIC$a$ [\citet{LF09}] with $p_{\lambda}(|\beta_j|)=(\lambda+1)|\beta
_j|/(\lambda+|\beta_j|)$ and MCP$+$ [\citet{Z09}] with all regularization
parameters being computed with 5-fold sparse generalized cross
validation; see Section \ref{subsec73} and Table \ref{tab4} therein
for detailed discussion on the choice of cross validation statistics.

The results of the simulations are summarized into three tables (see
Table~\ref{tab2} in the main text and Tables 2 and
3 in the supplementary material [\citet{BFJ11}]) where we reported the median
prediction error (PE)
%
\begin{equation}
\mathbb{P}_n[\exp\{-\bbeta^{*T} \bX\} - \exp\{ -\hat\bbeta{}^T
\bX\}]^2,
\end{equation}
where $\mathbb{P}_n$ stands for the empirical probability measure. We
also report the median number of nonzero parameters estimated in the
set $\mathcal{M}_*$ as the number of true positives TP. Furthermore, we
summarize the median number of nonzero estimates of the set $\mathcal
{M}_*^c$ as the number of false positives FP.

\begin{table}
\tabcolsep=0pt
\caption{Simulation results for $p \geq n$ under correlation
settings ranging from 0.25 to 0.90 with medium prediction error (MPE),
\# of true positives (TP), \# of false positives (FP) and standard
deviation in parenthesis of each estimate}
\label{tab2}
{\fontsize{8.8pt}{10.6pt}\selectfont{
\begin{tabular*}{\tablewidth}{@{\extracolsep{\fill}}ld{1.11}k{1.7}k{2.7}@{\hspace*{12pt}}d{1.11}k{1.7}k{2.8}@{}}
\hline
&\multicolumn{1}{c}{\textbf{MPE}}&\multicolumn{1}{c}{\textbf{TP}}
&\multicolumn{1}{c}{\textbf{FP}\hspace*{12pt}}&\multicolumn{1}{c}{\hspace*{-5pt}\textbf{MPE}}
&\multicolumn{1}{c}{\textbf{TP}}&\multicolumn{1}{c@{}}{\textbf{FP}}\\
\hline
& \multicolumn{6}{c@{}}{Settings of $n=100,p=100, s=4$}\\
& \multicolumn{3}{c}{Case $\rho=0.25$} &
\multicolumn{3}{c@{}}{Case $\rho=0.5$}\\
[2pt]
Oracle &0.0154\mbox{ (1.27\tabnoteref{ta})}&4&0&0.0215\mbox{ (1.97\tabnoteref{ta})} &4&0\\
LASSO &0.0178\mbox{ (1.26)}&4,\mbox{ (1.61)}&2,\mbox{ (33.34)} & 0.0284\mbox{ (2.12)}&4,\mbox{ (1.52)}&13,\mbox{ (33.02)}
\\
SCAD &0.0161\mbox{ (1.24)}&4,\mbox{ (1.61)}&2,\mbox{ (34.21)}& 0.0223\mbox{ (2.03)}&4,\mbox{ (1.52)}&13,\mbox{ (35.56)}\\
SICa &0.0190\mbox{ (1.27)}&3,\mbox{ (1.48)}&2,\mbox{ (26.11)}& 0.0275\mbox{ (2.43)}&3,\mbox{ (1.44)}&9,\mbox{ (21.54)}\\
MCP$+$&0.0166\mbox{ (1.22)} &3,\mbox{ (1.71)}&2,\mbox{ (32.49)}&0.0271\mbox{ (2.33)} &4,\mbox{ (1.54)}&24,\mbox{ (34.62)}\\
[4pt]
& \multicolumn{3}{c}{Case $\rho=0.75$} & \multicolumn
{3}{c@{}}{Case $\rho=0.9$} \\[2pt]
Oracle &0.0322\mbox{ (2.05\tabnoteref{ta})}&4&0&0.0538\mbox{ (4.43\tabnoteref{ta})}&4&0\\
LASSO &0.0371\mbox{ (2.42)}&3,\mbox{ (1.14)}&12,\mbox{ (31.21)} & 0.0665\mbox{ (4.62)}&2,\mbox{ (1.48)}&13,\mbox{ (32.16)}
\\
SCAD &0.0326\mbox{ (2.12)}&4,\mbox{ (1.14)}&12,\mbox{ (31.53)}& 0.0549\mbox{ (3.36)}&3,.5\mbox{ (1.49)}&8,\mbox{ (31.11)}\\
SICa &0.0343\mbox{ (2.27)}&2,\mbox{ (1.30)}&3,\mbox{ (18.41)}& 0.0566\mbox{ (3.26)}&2,\mbox{ (1.32)}&6,\mbox{ (24.42)}\\
MCP$+$&0.0326\mbox{ (2.21)} &3,.5\mbox{ (1.22)}&12,\mbox{ (32.31)}&0.0558\mbox{ (3.44)}
&2,.5\mbox{ (1.29)}&15,\mbox{ (29.68)}\\
[6pt]
& \multicolumn{6}{c@{}}{Settings of $n=100,p=1\mbox{,}000, s=4$}\\
& \multicolumn{3}{c}{Case $\rho=0.25$} &
\multicolumn{3}{c@{}}{Case $\rho=0.5$}\\[2pt]
Oracle &0.0154\mbox{ (1.27\tabnoteref{ta})}&4&0 & 0.0215\mbox{ (1.97\tabnoteref{ta})}&4&0\\
LASSO &0.0201\mbox{ (1.38)}&4,\mbox{ (0.85)}&23,\mbox{ (371.8)} &
0.0383\mbox{ (3.16)}&3,.5\mbox{ (1.23)}&45,\mbox{ (532.1)} \\
SCAD &0.0162\mbox{ (1.25)}&4,\mbox{ (0.83)}&15,\mbox{ (323.4)}& 0.0281\mbox{ (2.12)}&4,\mbox{ (1.12)}&36,\mbox{ (430.3)}\\
SICa &0.0189\mbox{ (1.17)}&3,.5\mbox{ (0.54)}&9,\mbox{ (120.5)}& 0.0492\mbox{ (3.18)}&3,\mbox{ (1.43)}&15,\mbox{ (319.4)}\\
MCP$+$&0.0192\mbox{ (1.23)} &4,\mbox{ (0.83)}&17,\mbox{ (345.5)}&0.0281\mbox{ (2.15)} &4,\mbox{ (1.12)}&36,\mbox{ (409.2)}\\[4pt]
& \multicolumn{3}{c}{Case $\rho=0.75$}
& \multicolumn{3}{c@{}}{Case $\rho=0.9$}\\[2pt]
Oracle &0.0322\mbox{ (2.05\tabnoteref{ta})}&4&0 &0.0538\mbox{ (4.43\tabnoteref{ta})}&4&0\\
LASSO &0.0497\mbox{ (3.16)}&3,\mbox{ (0.44)}&96,\mbox{ (306.5)} & 0.0703\mbox{ (4.24)}&
3,\mbox{ (1.54)}&97,\mbox{ (411.5)} \\
SCAD &0.0358\mbox{ (2.45)}&4,\mbox{ (0.34)}&85,\mbox{ (250.7)}& 0.0583\mbox{ (4.13)}&4,\mbox{ (1.51)}&67,\mbox{ (380.9)}\\
SICa &0.0372\mbox{ (2.15)}&2,\mbox{ (1.30)}&90,.5\mbox{ (90.3)}& 0.0546\mbox{ (3.98)}&1,\mbox{ (1.78)}&30,\mbox{ (354.1)}\\
MCP$+$&0.0361\mbox{ (2.77)} &3,.5\mbox{ (1.14)}&90,\mbox{ (320.4)}&0.0592\mbox{ (4.25)}
&3,.5\mbox{ (1.58)}&98,\mbox{ (402.3)}\\
[6pt]
& \multicolumn{6}{c@{}}{Settings of $n=100,p=5\mbox{,}000, s=4$}\\
& \multicolumn{3}{c}{Case $\rho=0.25$}
&\multicolumn{3}{c@{}}{Case $\rho=0.5$}\\[2pt]
Oracle &0.0154\mbox{ (1.27\tabnoteref{ta})}&4&0 &0.0215\mbox{ (1.97\tabnoteref{ta})}&4&0\\
LASSO &0.0220\mbox{ (1.49)}&4,\mbox{ (1.05)}&68,\mbox{ (398.1)} &
0.0462\mbox{ (4.05)}&3,.5\mbox{ (1.64)}&33,\mbox{ (206.8)} \\
SCAD &0.0170\mbox{ (1.28)}&4,\mbox{ (1.05)}&67,\mbox{ (298.2)}&
0.0328\mbox{ (3.15)}&3,.5\mbox{ (1.56)}&21,.5\mbox{ (205.4)}\\
SICa &0.0195\mbox{ (1.19)}&2,.5\mbox{ (1.17)}&14,\mbox{ (345.7)}& 0.0285\mbox{ (3.35)}&4,\mbox{ (1.41)}&30,\mbox{ (323.3)}\\
MCP$+$&0.0188\mbox{ (1.29)} &3,\mbox{ (1.10)}&67,\mbox{ (298.2)}&0.0358\mbox{ (2.85)}
&3,.5\mbox{ (1.51)}&73,.5\mbox{ (348.7)}\\[4pt]
& \multicolumn{3}{c}{Case $\rho=0.75$}
&\multicolumn{3}{c@{}}{Case $\rho=0.9$} \\[2pt]
Oracle &0.0322\mbox{ (2.05\tabnoteref{ta})}&4&0 &0.0538\mbox{ (4.43\tabnoteref{ta})}&4&0\\
LASSO &0.0567\mbox{ (5.02)}&3,\mbox{ (1.73)}&23,\mbox{ (250.5)} & 0.0865\mbox{ (4.52)}&
2,\mbox{ (1.23)}&59,\mbox{ (208.8)} \\
SCAD &0.0360\mbox{ (2.31)}&4,\mbox{ (1.51)}&18,\mbox{ (234.7)}& 0.0596\mbox{ (4.12)}&4,\mbox{ (0.89)}&49,\mbox{ (105.4)}\\
SICa &0.0385\mbox{ (2.13)}&2,.5\mbox{ (1.30)}&3,\mbox{ (225.2)}& 0.0602\mbox{ (4.92)}&3,\mbox{ (0.45)}&46,\mbox{ (90.3)}\\
MCP$+$&0.0392\mbox{ (2.82)} &4,\mbox{ (1.74)}&4,\mbox{ (326.2)}&0.0578\mbox{ (4.33)} &4,\mbox{ (0.89)}&11,\mbox{ (217.1)}\\
\hline
\end{tabular*}
\tabnotetext[\mbox{$\sharp$}]{ta}{stands for column of standard
deviation${}\times{}$100.}}}%
\end{table}

Table \ref{tab2} summarizes three $p\geq n$ examples, where especially
the last two
stress the strengths of the methods when $p \gg n$ and spectra of the
design matrix is high; see Table 1 in the supplementary
material [\citet{BFJ11}]. All
four methods work quite well, where LASSO has higher PE than the rest,
with SCAD and MCP$+$ performing quite closely to each other. SICa
performs worse than others, always loosing a number of TPs. The case of
$\rho=0.90$ affects all methods in bigger prediction error and smaller
number of TP, where the jump is the largest in LASSO penalty. SCAD and
MCP keep their performance similarly to the oracle one through all
examples, hence verifying the strengths of nonconvex penalties. For
more detailed discussions and results when the oracle estimator fails,
when the censoring rate is too high and assessing the relative
estimation efficiency of LASSO estimator with respect to SCAD, SICa and
MCP$+$, we direct you to the supplementary material [\citet{BFJ11}]
for this paper.



\subsection{Real data example} \label{subsec73}

To demonstrate the strength of the proposed methodology, in this
section, we present gene association study with respect to the survival
time of non-Hodgkin's lymphoma.
Genetic mechanisms responsible for the clinical heterogeneity of
follicular lymphoma are still unknown. \citet{D04} have collected gene
expression data on 191 biopsy specimens obtained from patients with
untreated follicular lymphoma. RNA was extracted from fresh-frozen
tumor-biopsy specimens and survival times, from 191 patients, who had
received a diagnosis between 1974 and 2001, which were obtained from
seven institutions and examined for gene expression with the use of
Affymetrix U133A and U133B microarrays. The median age at diagnosis was
51 years (range, 23 to 81), and the median follow-up time was 6.6 years
(range, less than 1.0 to 28.2).
The dataset was obtained from \url{http://llmpp.nih.gov/FL}.

The full cohort study included 44,187 probe expressions values out
which only 34,188 were properly annotated. Among these, some received
multiple (2--7) measurements per gene. We took the median value as a
unique representative and were left with 17,118 different genes
presented. We separated the dataset into training and testing sets with
80\% and 20\% of censored samples, respectively. The censoring rate of
50\% was kept in each of the training and testing samples. Recorded for
each individual are follow up time, indicator of the status at the
follow up time and measurements of expression value for each Affymetrix
probe set.

The classical $L$ fold cross-validation is defined as
\[
\operatorname{CV}(\lambda)=\sum_{k=1}^{L} \bigl\{ l\bigl(\hat\bbeta{}^{(-k)}_{\lambda}\bigr)
-l^{(-k)}\bigl(\hat\bbeta{^{(-k)}}_{\lambda}\bigr)\bigr\},
\]
where $l$ stands for the partial likelihood and $l^{(-k)}$ for the
partial likelihood evaluated\vadjust{\goodbreak} without the $k$th subset and similarly
$\hat\bbeta{}^{(-k)}_\lambda$ for the penalized estimator derived without
using the $k$th subset.
The measure of information contained in the full Cox partial likelihood
is biased with respect to the number of nonzero elements and proper
normalization is needed. The method of generalized cross validation
proposed by \citet{FL02} works very well for small $p$ but fails for
large $p$ because of its dependence on the inverse of the Hessian
matrix of the partial likelihood. This inspired us to define a~sparse
approximation to the generalized cross-validation
as
\[
\operatorname{SGCV}(\lambda)=
\sum_{k=1}^{L} \biggl( \frac{l(\hat\bbeta{}^{(-k)}_{\lambda})}{n\{
1-\hat
s_\lambda/n\}^2} - \frac{l^{(-k)}(\hat\bbeta{}^{(-k)}_{\lambda
})}{n^{(-k)}\{1-\hat s_\lambda/n^{(-k)}\}^2}\biggr),
\]
where $\hat s_{\lambda}=\|\hat\beta_{\lambda}^{(-k)}\|_0$ and
$n^{(-k)}$ stands for the sample size of the whole set without the
$k$th subset. Then, we choose the regularization parameter as
\[
\hat{\lambda}=\mathop{\arg\min}_{\lambda\dvtx \hat s_{\lambda} < n}
\operatorname{SGCV}(\lambda).
\]

We applied 5-fold cross validation on the test set and evaluated its
performance on the training set. The Nelson--Aalen estimate of the
cumulative hazard rate function was used. The results\vadjust{\goodbreak} are summarized in
Table \ref{tab4} and show a big difference between the classical CV
statistics and generalized one. The CV, being not scaled to the number
of nonzero elements always prefers models with bigger number of
nonzeros. Note that $\hat s >n$, for small~$\lambda$, is caused by the
artifact of ICA algorithm.

\begin{table}
\caption{Data summary with number of nonzero elements reported on
the whole data set and prediction error and its standard
deviation${}\times{}$100 comparisons reported on the training
set [Dave et~al. (\citeyear{D04})]}\label{tab4}
\begin{tabular*}{\tablewidth}{@{\extracolsep{\fill}}lcccc@{}}
\hline
&\textbf{LASSO}&\textbf{SCAD}&\textbf{SICA}&\textbf{MCP}$\bolds{+}$\\
\hline
CV & &&&\\
\# of nonzeros & 2145 &653&0&154\\
Prediction error & 0.1516 (1.51) &0.1276 (1.60)&0.1898 (--)&0.1743 (1.45)\\
[4pt]
SGCV&&&&\\
\# of nonzeros &24 &26&0&13\\
Prediction error & 0.0812 (1.03) &0.0643 (1.02)&0.1898 (--)&0.1043 (0.78)\\
\hline
\end{tabular*}
\end{table}

The SICa penalty completely fails in this example. It detects nonzeros
only in 3 grid points with the number of nonzeros as $2,3$ and $879$.
Both CV methods fail to pick up the optimal one among the three points
and choose the fourth one, which lead to no signal detection. This is
not unexpected, since in all simulations SICa was always picking the
least number of TP$+$FP; see Table \ref{tab2}.

\begin{table}
\caption{Data estimation summary of the genes selected by the
sparse generalized cross validation with standard deviation$\times
100^\ddag$ reported in the parenthesis}\label{tab5}
\begin{tabular*}{\tablewidth}{@{\extracolsep{\fill}}lk{2.14}k{2.14}k{2.14}}
\hline
\textbf{Gene annotation}&\multicolumn{1}{c}{\textbf{LASSO}}
&\multicolumn{1}{c}{\textbf{SCAD}}&\multicolumn{1}{c@{}}{\textbf{MCP}$\bolds{+}$}\\
\hline
FOSB\ (BC036724)&-0.,0093\ (2.34^\ddag)&\multicolumn{1}{c}{$\times$} &-0.,0027\ (1.54^\ddag)\\
GABRA6\ (AK090735)& 0.,0070\ (0.56)^{*}&0.,0150\ (1.00^\ddag)^{***}&\multicolumn{1}{c@{}}{$\times$}\\
GHRH\ (AW$\_{}134884$)&\multicolumn{1}{c}{$\times$} &-0.,0489 \ (1.39)^{**}&\multicolumn{1}{c@{}}{$\times$}\\
GNGT\_1\ (BC030956)& \multicolumn{1}{c}{$\times$}&-0.,0041 \ (0.46)&\multicolumn{1}{c@{}}{$\times$}\\
HIST1H1E\ (BU603483)&-0.,0026\ (1.98) &-0.,0032 \ (1.41)&\multicolumn{1}{c@{}}{$\times$}\\
HIST1H2AE\ (BE741093)&\multicolumn{1}{c}{$\times$} &-0.,0137 \ (0.41)^{**}&\multicolumn{1}{c@{}}{$\times$}\\
IFNA2\ (NM$\_{}000605$)&\multicolumn{1}{c}{$\times$} &0.,0095 \ (1.29)&\multicolumn{1}{c@{}}{$\times$}\\
IMPG1\ (NM$\_{}001563$)& -0.,0168\ (2.56)&-0.,0116 \ (0.81)&\multicolumn{1}{c@{}}{$\times$}\\
MATN3\ (NM$\_{}002381$)&0.,0206\ (0.89)^{**}&0.,0301
\ (0.36)^{***}&0.,0065\ (1.25)\\
RTH\ (NM$\_{}000315$)&-0.,0032\ (0.85) &-0.,0177 \ (0.74)^{**}&\multicolumn{1}{c@{}}{$\times$}\\
RAG2\ (NM$\_{}000536$)& -8.,8728\mbox{e--05}\ (1.56)&\multicolumn{1}{c}{$\times$}&\multicolumn{1}{c@{}}{$\times$}\\
SCN9A\ (NM$\_{}002977$)&0.,0049\ (0.87)&\multicolumn{1}{c}{$\times$}& 8.,5785\mbox{e--04}\ (1.56)\\
CXCL5\ (NM$\_{}002994$)& \multicolumn{1}{c}{$\times$}& 0.,0026 \ (1.69)&\multicolumn{1}{c@{}}{$\times$}\\
SH3BGR\ (BM725357)&-0.,0125\ (0.25)^{***} &\multicolumn{1}{c}{$\times$}&\multicolumn{1}{c@{}}{$\times$}\\
HIST1H3B\ (NM$\_{}006770$)&\multicolumn{1}{c}{$\times$} &-0.,0029 \ (0.81)&\multicolumn{1}{c@{}}{$\times$}\\
MARCO\ (BP872375)&0.,0013 \ (2.54)&\multicolumn{1}{c}{$\times$}&\multicolumn{1}{c@{}}{$\times$}\\
CLCA3\ (NM$\_{}004921$) &0.,0172\ (0.85)^{*}&0.,0170 \ (0.71)^{*}&
0.,0171\ (0.54)^{**}\\
SEMA3A\ (XM$\_{}376647$)& -0.,0049\ (1.15)&\multicolumn{1}{c}{$\times$}& -3.,8781\mbox{e--05}\ (0.76)\\
KIAA0861\ (BX694003)&-0.,0261\ (0.96)^{*}&-0.,0181 \ (0.74)^{**}&
-0.,0170\ (0.58)^{**}\\
FSCN2\ (NM$\_{}012418$)&0.,0136\ (1.25)& 0.,0194 \ (1.44)& 0.,0058\ (1.12)\\
DKFZP566K0\ (ALO50040)& -0.,0025\ (1.36)&\multicolumn{1}{c}{$\times$}&\multicolumn{1}{c@{}}{$\times$}\\
MORC\ (BC050307)&0.,0204\ (0.75)^{*}&0.,0204 \ (1.00)^{*}& 0.,0165\ (0.94)\\
C14orf105\ (ALO1512)& 0.,0021\ (1.47)&\multicolumn{1}{c}{$\times$}&\multicolumn{1}{c@{}}{$\times$}\\
SAGE1\ (NM$\_{}018667$)& \multicolumn{1}{c}{$\times$}&0.,0012 \ (2.56)&\multicolumn{1}{c@{}}{$\times$}\\
C6orf103\ (AL832192)& \multicolumn{1}{c}{$\times$}& 0.,0023 \ (1.20)&\multicolumn{1}{c@{}}{$\times$}\\
FLJ13841\ (AK023903)&0.,0146\ (0.35)^{**}&0.,0146 \ (0.49)^{**}&
0.,0129\ (0.47)^{*}\\
FLJ22655\ (BC042888)& \multicolumn{1}{c}{$\times$}&0.,0028 \ (0.64)&\multicolumn{1}{c@{}}{$\times$}\\
FLJ21934\ (AY358727)&-0.,0127\ (0.55)^{*}&-0.,0125 \ (0.49)^{**}&
-0.,0079\ (0.32)^{*}\\
KIAA1912\ (AB067499)&\multicolumn{1}{c}{$\times$} &-0.,0013 \ (0.65)&\multicolumn{1}{c@{}}{$\times$}\\
FLJ40298\ (NM\_173486)& 0.,0307\ (0.98)^{***}&0.,0316 \ (1.20)^{**}&\multicolumn{1}{c@{}}{$\times
$}\\
MGC33951\ (BC029537)&\multicolumn{1}{c}{$\times$} &-0.,0042 \ (1.44)&\multicolumn{1}{c@{}}{$\times$}\\
NALP4\ (AF479747) &-0.,0059\ (0.56)&-0.,0059 \ (0.76)& -0.,0062\ (0.23)^{*}\\
FLJ46154\ (AK128035)&0.,0185\ (0.89)^{*}&0.,0185 \ (0.94)^{*}&
0.,0141\ (1.58)\\
MGC50372\ (BX647272)&5.,3676\mbox{e--04} \ (1.68)&\multicolumn{1}{c}{$\times$}&\multicolumn{1}{c@{}}{$\times$}\\
LOC285016\ (XM$\_{}211736$) &0.,0182\ (2.56)&0.,0182 \ (2.25)& 0.,0147\ (2.15)\\
\hline
\end{tabular*}
\legend{Superscripts $^{***}$, $^{**}$, $^*$ are decodings of significance
values.}
\end{table}

Table \ref{tab5} depicts the estimation results of the sparse
generalized cross validation method with LASSO, SCAD and MCP$+$
penalties. All three penalties yield sign consistency of estimated
coefficients among the selected gene sets. Note that the relative
rankings of estimated corresponding coefficients are different among
all methods. For example, gene FLJ40298 has the biggest absolute size
in the SCAD penalty, it is ranked number 5 among those coefficients
produced by LASSO penalty and it is not even selected in the MCP$+$
penalty. Interestingly, the common set of genes selected by LASSO and
SCAD has very consistent estimated coefficients.
For most genes MCP results in smaller estimated values than SCAD and LASSO.

\section{Discussion}

We have studied penalized log partial likelihood methods for ultra-high
dimensional variable selection for Cox's regression models. With
nonconcave penalties, we have shown
that such methods have model selection consistency with oracle properties
even for NP-dimensionality. We have established that oracle properties
hold with probability converging to one exponentially fast, and that
the rate explicitly depends on the real and intrinsic dimensionality
$p$ and~$s$, respectively. We have also developed an exponential
inequality for deviations of a counting process from its compensator.
Results for LASSO penalty were obtained as a special case. It confirms
explicitly that folded concave penalties allow for far weaker
correlation structure than LASSO penalty.
Furthermore, the asymptotic normality was proved, results of which can
be used to construct confidence intervals of the estimated
coefficients.\vadjust{\goodbreak}

%

%

%


\begin{supplement}
\stitle{Supplementary material for ``Regularization
for Cox's proportional hazards model with
NP-dimensionality''}
\slink[doi]{10.1214/11-AOS911SUPP} 
\sdatatype{.pdf}
\sfilename{aos911\_supp.pdf}
\sdescription{In the Supplementary Material [\citet{BFJ11}] we give
additional results of our simulation study, we specify the statements
and  detailed proofs of  technical Lemmas
2.1--2.3 and give complete proofs of
Theorems~\ref{theo21}, \ref{theo411},
\mbox{\ref{cor2}--\ref{theo61}}. We present the details of the ICA
algorithm of the Section~\ref{sec7} together with new simulation
settings were we  increased the censoring rate and/or increased the
number of significant variables $s$, and with discussion on the
relative estimation efficiency of the penalized methods. We develop
results on the growth of the $L_2$ norm of the score vector~%
$U_n(\bbeta_1^*)$ and of the matrix $\int_0^\tau \bV(\bbeta_1^*,t) \,d
\bar M(t)$. Moreover we establish a result on the asymptotic behavior
of vector $\hat \bbeta_1^*$ when $s=o(n^{1/3})$ diverging with $n$. The
main tools used are  the theory of martingales [\citet{FH05}] and the
results of various matrix norms of  Lemmas \ref{le92a},~\ref{le93} and~2.1--2.3.}
\end{supplement}

%

\printaddresses


\begin{thebibliography}{40}

\bibitem[\protect\citeauthoryear{Andersen and Gill}{1982}]{AG82}
%
\begin{barticle}[mr]
\bauthor{\bsnm{Andersen},~\bfnm{P.~K.}\binits{P.~K.}} \AND
\bauthor{\bsnm{Gill},~\bfnm{R.~D.}\binits{R.~D.}}
(\byear{1982}).
\btitle{Cox's regression model for counting processes: A large sample study}.
\bjournal{Ann. Statist.}
\bvolume{10}
\bpages{1100--1120}.
\bid{issn={0090-5364}, mr={0673646}}
\bptok{imsref}%
\end{barticle}
%
\endbibitem

\bibitem[\protect\citeauthoryear{Bertsekas}{2003}]{B03}
%
\begin{bmisc}[mr]
\bauthor{\bsnm{Bertsekas},~\bfnm{Dimitri~P.}\binits{D.~P.}}
(\byear{2003}).
\bhowpublished{Nonlinear programming. Athena Scientific,
Nashua, NH.}
\bptok{imsref}%
\end{bmisc}
%
\endbibitem

\bibitem[\protect\citeauthoryear{Bhatia}{1997}]{B97}
%
\begin{bbook}[mr]
\bauthor{\bsnm{Bhatia},~\bfnm{Rajendra}\binits{R.}}
(\byear{1997}).
\btitle{Matrix Analysis}.
\bseries{Graduate Texts in Mathematics}
\bvolume{169}.
\bpublisher{Springer}, \baddress{New York}.
\bid{mr={1477662}}
\bptok{imsref}%
\end{bbook}
%
\endbibitem

\bibitem[\protect\citeauthoryear{Bickel, Ritov and Tsybakov}{2009}]{BRT09}
%
\begin{barticle}[mr]
\bauthor{\bsnm{Bickel},~\bfnm{Peter~J.}\binits{P.~J.}},
\bauthor{\bsnm{Ritov},~\bfnm{Ya'acov}\binits{Y.}} \AND
\bauthor{\bsnm{Tsybakov},~\bfnm{Alexandre~B.}\binits{A.~B.}}
(\byear{2009}).
\btitle{Simultaneous analysis of lasso and {D}antzig selector}.
\bjournal{Ann. Statist.}
\bvolume{37}
\bpages{1705--1732}.
\bid{doi={10.1214/08-AOS620}, issn={0090-5364}, mr={2533469}}
\bptok{imsref}%
\end{barticle}
%
\endbibitem

\bibitem[\protect\citeauthoryear{Bradic, Fan and Wang}{2011}]{BFW09}
%
\begin{barticle}[mr]
\bauthor{\bsnm{Bradic},~\bfnm{Jelena}\binits{J.}},
\bauthor{\bsnm{Fan},~\bfnm{Jianqing}\binits{J.}} \AND
\bauthor{\bsnm{Wang},~\bfnm{Weiwei}\binits{W.}}
(\byear{2011}).
\btitle{Penalized composite quasi-likelihood for ultrahigh dimensional variable
selection}.
\bjournal{J. R. Stat. Soc. Ser. B Stat. Methodol.}
\bvolume{73}
\bpages{325--349}.
\bid{doi={10.1111/j.1467-9868.2010.00764.x}, issn={1369-7412}, mr={2815779}}
\bptok{imsref}%
\end{barticle}
%
\endbibitem

\bibitem[\protect\citeauthoryear{Bradic, Fan and Jiang}{2011}]{BFJ11}
%
\begin{bmisc}[auto:STB|2011/11/09|09:54:39]
\bauthor{\bsnm{Bradic},~\bfnm{J.}\binits{J.}},
\bauthor{\bsnm{Fan},~\bfnm{J.}\binits{J.}} \AND
\bauthor{\bsnm{Jiang},~\bfnm{J.}\binits{J.}}
(\byear{2011}).
\bhowpublished{Supplement to ``Regularization for Cox's proportional
hazards model with NP-dimensionality.''
\href{http://dx.doi.org/10.1214/11-AOS911SUPP}{DOI:10.1214/11-AOS911SUPP}.}
\bptok{imsref}%
\end{bmisc}
%
\endbibitem

\bibitem[\protect\citeauthoryear{Bunea, Tsybakov and Wegkamp}{2007}]{BTW07}
%
\begin{barticle}[mr]
\bauthor{\bsnm{Bunea},~\bfnm{Florentina}\binits{F.}},
\bauthor{\bsnm{Tsybakov},~\bfnm{Alexandre}\binits{A.}} \AND
\bauthor{\bsnm{Wegkamp},~\bfnm{Marten}\binits{M.}}
(\byear{2007}).
\btitle{Sparsity oracle inequalities for the {L}asso}.
\bjournal{Electron. J. Stat.}
\bvolume{1}
\bpages{169--194 (electronic)}.
\bid{doi={10.1214/07-EJS008}, issn={1935-7524}, mr={2312149}}
\bptok{imsref}%
\end{barticle}
%
\endbibitem

\bibitem[\protect\citeauthoryear{Cai et~al.}{2005}]{CFLZ05}
%
\begin{barticle}[mr]
\bauthor{\bsnm{Cai},~\bfnm{Jianwen}\binits{J.}},
\bauthor{\bsnm{Fan},~\bfnm{Jianqing}\binits{J.}},
\bauthor{\bsnm{Li},~\bfnm{Runze}\binits{R.}} \AND
\bauthor{\bsnm{Zhou},~\bfnm{Haibo}\binits{H.}}
(\byear{2005}).
\btitle{Variable selection for multivariate failure time data}.
\bjournal{Biometrika}
\bvolume{92}
\bpages{303--316}.
\bid{doi={10.1093/biomet/92.2.303}, issn={0006-3444}, mr={2201361}}
\bptok{imsref}%
\end{barticle}
%
\endbibitem

\bibitem[\protect\citeauthoryear{Candes and Tao}{2007}]{CT07}
%
\begin{barticle}[mr]
\bauthor{\bsnm{Candes},~\bfnm{Emmanuel}\binits{E.}} \AND
\bauthor{\bsnm{Tao},~\bfnm{Terence}\binits{T.}}
(\byear{2007}).
\btitle{The {D}antzig selector: Statistical estimation when {$p$} is much
larger than {$n$}}.
\bjournal{Ann. Statist.}
\bvolume{35}
\bpages{2313--2351}.
\bid{doi={10.1214/009053606000001523}, issn={0090-5364}, mr={2382644}}
\bptok{imsref}%
\end{barticle}
%
\endbibitem

\bibitem[\protect\citeauthoryear{Daubechies, Defrise and
De~Mol}{2004}]{DDD04}
%
\begin{barticle}[mr]
\bauthor{\bsnm{Daubechies},~\bfnm{Ingrid}\binits{I.}},
\bauthor{\bsnm{Defrise},~\bfnm{Michel}\binits{M.}} \AND
\bauthor{\bsnm{De~Mol},~\bfnm{Christine}\binits{C.}}
(\byear{2004}).
\btitle{An iterative thresholding algorithm for linear inverse
problems with a
sparsity constraint}.
\bjournal{Comm. Pure Appl. Math.}
\bvolume{57}
\bpages{1413--1457}.
\bid{doi={10.1002/cpa.20042}, issn={0010-3640}, mr={2077704}}
\bptok{imsref}%
\end{barticle}
%
\endbibitem

\bibitem[\protect\citeauthoryear{Dave et~al.}{2004}]{D04}
%
\begin{barticle}[pbm]
\bauthor{\bsnm{Dave},~\bfnm{Sandeep~S.}\binits{S.~S.}~\bsuffix{et al.}}
(\byear{2004}).
\btitle{Prediction of survival in follicular lymphoma based on molecular
features of tumor-infiltrating immune cells}.
\bjournal{N. Engl. J. Med.}
\bvolume{351}
\bpages{2159--2169}.
\bid{doi={10.1056/NEJMoa041869}, issn={1533-4406}, pii={351/21/2159},
pmid={15548776}}
\bptok{imsref}%
\end{barticle}
%
\endbibitem

\bibitem[\protect\citeauthoryear{de~la Pe{\~n}a}{1999}]{dlP99}
%
\begin{barticle}[mr]
\bauthor{\bparticle{de~la} \bsnm{Pe{\~n}a},~\bfnm{Victor~H.}\binits{V.~H.}}
(\byear{1999}).
\btitle{A general class of exponential inequalities for martingales and
ratios}.
\bjournal{Ann. Probab.}
\bvolume{27}
\bpages{537--564}.
\bid{doi={10.1214/aop/1022677271}, issn={0091-1798}, mr={1681153}}
\bptok{imsref}%
\end{barticle}
%
\endbibitem

\bibitem[\protect\citeauthoryear{Du, Ma and Liang}{2010}]{DML10}
%
\begin{barticle}[mr]
\bauthor{\bsnm{Du},~\bfnm{Pang}\binits{P.}},
\bauthor{\bsnm{Ma},~\bfnm{Shuangge}\binits{S.}} \AND
\bauthor{\bsnm{Liang},~\bfnm{Hua}\binits{H.}}
(\byear{2010}).
\btitle{Penalized variable selection procedure for {C}ox models with
semiparametric relative risk}.
\bjournal{Ann. Statist.}
\bvolume{38}
\bpages{2092--2117}.
\bid{doi={10.1214/09-AOS780}, issn={0090-5364}, mr={2676884}}
\bptnote{check year}%
\bptok{imsref}%
\end{barticle}
%
\endbibitem

\bibitem[\protect\citeauthoryear{Fan and Li}{2001}]{FL01}
%
\begin{barticle}[mr]
\bauthor{\bsnm{Fan},~\bfnm{Jianqing}\binits{J.}} \AND
\bauthor{\bsnm{Li},~\bfnm{Runze}\binits{R.}}
(\byear{2001}).
\btitle{Variable selection via nonconcave penalized likelihood and its oracle
properties}.
\bjournal{J. Amer. Statist. Assoc.}
\bvolume{96}
\bpages{1348--1360}.
\bid{doi={10.1198/016214501753382273}, issn={0162-1459}, mr={1946581}}
\bptok{imsref}%
\end{barticle}
%
\endbibitem

\bibitem[\protect\citeauthoryear{Fan and Li}{2002}]{FL02}
%
\begin{barticle}[mr]
\bauthor{\bsnm{Fan},~\bfnm{Jianqing}\binits{J.}} \AND
\bauthor{\bsnm{Li},~\bfnm{Runze}\binits{R.}}
(\byear{2002}).
\btitle{Variable selection for {C}ox's proportional hazards model and frailty
model}.
\bjournal{Ann. Statist.}
\bvolume{30}
\bpages{74--99}.
\bid{doi={10.1214/aos/1015362185}, issn={0090-5364}, mr={1892656}}
\bptok{imsref}%
\end{barticle}
%
\endbibitem

\bibitem[\protect\citeauthoryear{Fan and Lv}{2011}]{FL10}
%
\begin{barticle}[auto:STB|2011/11/09|09:54:39]
\bauthor{\bsnm{Fan},~\bfnm{J.}\binits{J.}} \AND
\bauthor{\bsnm{Lv},~\bfnm{J.}\binits{J.}}
(\byear{2011}).
\btitle{Non-concave penalized likelihood with NP-dimensionality}.
\bjournal{IEEE Trans. Inform. Theory}
\bvolume{57}
\bpages{5467--5484}.
\bptok{imsref}%
\end{barticle}
%
\endbibitem

\bibitem[\protect\citeauthoryear{Fleming and Harrington}{1991}]{FH05}
%
\begin{bbook}[mr]
\bauthor{\bsnm{Fleming},~\bfnm{Thomas~R.}\binits{T.~R.}} \AND
\bauthor{\bsnm{Harrington},~\bfnm{David~P.}\binits{D.~P.}}
(\byear{1991}).
\btitle{Counting Processes and Survival Analysis}.
\bpublisher{Wiley}, \baddress{New York}.
\bid{mr={1100924}}
\bptnote{check year}%
\bptok{imsref}%
\end{bbook}
%
\endbibitem


\bibitem[\protect\citeauthoryear{Friedman, Hastie and
Tibshirani}{2010}]{FHT10}
%
\begin{barticle}[pbm]
\bauthor{\bsnm{Friedman},~\bfnm{Jerome}\binits{J.}},
\bauthor{\bsnm{Hastie},~\bfnm{Trevor}\binits{T.}} \AND
\bauthor{\bsnm{Tibshirani},~\bfnm{Rob}\binits{R.}}
(\byear{2010}).
\btitle{Regularization paths for generalized linear models via coordinate
descent}.
\bjournal{Journal of Statistical Software}
\bvolume{33}
\bpages{1--22}.
\bid{issn={1548-7660}, mid={NIHMS201118}, pmcid={2929880}, pmid={20808728}}
\bptok{imsref}%
\end{barticle}
%
\endbibitem

\bibitem[\protect\citeauthoryear{Friedman et~al.}{2007}]{FHHT07}
%
\begin{barticle}[mr]
\bauthor{\bsnm{Friedman},~\bfnm{Jerome}\binits{J.}},
\bauthor{\bsnm{Hastie},~\bfnm{Trevor}\binits{T.}},
\bauthor{\bsnm{H{\"o}fling},~\bfnm{Holger}\binits{H.}} \AND
\bauthor{\bsnm{Tibshirani},~\bfnm{Robert}\binits{R.}}
(\byear{2007}).
\btitle{Pathwise coordinate optimization}.
\bjournal{Ann. Appl. Stat.}
\bvolume{1}
\bpages{302--332}.
\bid{doi={10.1214/07-AOAS131}, issn={1932-6157}, mr={2415737}}
\bptok{imsref}%
\end{barticle}
%
\endbibitem

\bibitem[\protect\citeauthoryear{Johnson}{2009}]{J09}
%
\begin{barticle}[mr]
\bauthor{\bsnm{Johnson},~\bfnm{Brent~A.}\binits{B.~A.}}
(\byear{2009}).
\btitle{On lasso for censored data}.
\bjournal{Electron. J. Stat.}
\bvolume{3}
\bpages{485--506}.
\bid{doi={10.1214/08-EJS322}, issn={1935-7524}, mr={2507457}}
\bptok{imsref}%
\end{barticle}
%
\endbibitem

\bibitem[\protect\citeauthoryear{Juditsky and Nemirovski}{2011}]{JN10}
%
\begin{bmisc}[auto:STB|2011/11/09|09:54:39]
\bauthor{\bsnm{Juditsky},~\bfnm{A.~B.}\binits{A.~B.}} \AND
\bauthor{\bsnm{Nemirovski},~\bfnm{A.~S.}\binits{A.~S.}}
(\byear{2011}).
\bhowpublished{Large deviations of vector-valued martingales in
2-smooth normed
spaces.
\textit{Ann. Appl. Probab.} To appear.
Available at \href{http://arxiv.org/abs/0809.0813}{arXiv:0809.0813}.}
\bptok{imsref}%
\end{bmisc}
%
\endbibitem

\bibitem[\protect\citeauthoryear{Kim, Choi and Oh}{2008}]{KCO08}
%
\begin{barticle}[auto:STB|2011/11/09|09:54:39]
\bauthor{\bsnm{Kim},~\bfnm{Y.}\binits{Y.}},
\bauthor{\bsnm{Choi},~\bfnm{H.}\binits{H.}} \AND
\bauthor{\bsnm{Oh},~\bfnm{H.}\binits{H.}}
(\byear{2008}).
\btitle{Smoothly clipped absolute deviation on high dimensions}.
\bjournal{J.~Amer. Statist. Assoc.}
\bvolume{103}
\bpages{1656--1673}.
\bptok{imsref}%
\end{barticle}
%
\endbibitem

\bibitem[\protect\citeauthoryear{Koltchinskii}{2009}]{K09}
%
\begin{barticle}[mr]
\bauthor{\bsnm{Koltchinskii},~\bfnm{Vladimir}\binits{V.}}
(\byear{2009}).
\btitle{The {D}antzig selector and sparsity oracle inequalities}.
\bjournal{Bernoulli}
\bvolume{15}
\bpages{799--828}.
\bid{doi={10.3150/09-BEJ187}, issn={1350-7265}, mr={2555200}}
\bptok{imsref}%
\end{barticle}
%
\endbibitem

\bibitem[\protect\citeauthoryear{Lv and Fan}{2009}]{LF09}
%
\begin{barticle}[mr]
\bauthor{\bsnm{Lv},~\bfnm{Jinchi}\binits{J.}} \AND
\bauthor{\bsnm{Fan},~\bfnm{Yingying}\binits{Y.}}
(\byear{2009}).
\btitle{A unified approach to model selection and sparse recovery using
regularized least squares}.
\bjournal{Ann. Statist.}
\bvolume{37}
\bpages{3498--3528}.
\bid{doi={10.1214/09-AOS683}, issn={0090-5364}, mr={2549567}}
\bptok{imsref}%
\end{barticle}
%
\endbibitem

\bibitem[\protect\citeauthoryear{Massart and Meynet}{2010}]{MM10}
%
\begin{bmisc}[auto:STB|2011/11/09|09:54:39]
\bauthor{\bsnm{Massart},~\bfnm{P.}\binits{P.}} \AND
\bauthor{\bsnm{Meynet},~\bfnm{C.}\binits{C.}}
(\byear{2010}).
\bhowpublished{An $l_1$ oracle inequality for the LASSO.
Available at \href{http://arxiv.org/abs/1007.4791}{arXiv:1007.4791}}.
\bptok{imsref}%
\end{bmisc}
%
\endbibitem

\bibitem[\protect\citeauthoryear{Meinshausen and B{\"u}hlmann}{2006}]{MB06}
%
\begin{barticle}[mr]
\bauthor{\bsnm{Meinshausen},~\bfnm{Nicolai}\binits{N.}} \AND
\bauthor{\bsnm{B{\"u}hlmann},~\bfnm{Peter}\binits{P.}}
(\byear{2006}).
\btitle{High-dimensional graphs and variable selection with the lasso}.
\bjournal{Ann. Statist.}
\bvolume{34}
\bpages{1436--1462}.
\bid{doi={10.1214/009053606000000281}, issn={0090-5364}, mr={2278363}}
\bptok{imsref}%
\end{barticle}
%
\endbibitem

\bibitem[\protect\citeauthoryear{Meinshausen and Yu}{2009}]{MY09}
%
\begin{barticle}[mr]
\bauthor{\bsnm{Meinshausen},~\bfnm{Nicolai}\binits{N.}} \AND
\bauthor{\bsnm{Yu},~\bfnm{Bin}\binits{B.}}
(\byear{2009}).
\btitle{Lasso-type recovery of sparse representations for high-dimensional
data}.
\bjournal{Ann. Statist.}
\bvolume{37}
\bpages{246--270}.
\bid{doi={10.1214/07-AOS582}, issn={0090-5364}, mr={2488351}}
\bptok{imsref}%
\end{barticle}
%
\endbibitem

\bibitem[\protect\citeauthoryear{Tibshirani}{1996}]{T96}
%
\begin{barticle}[mr]
\bauthor{\bsnm{Tibshirani},~\bfnm{Robert}\binits{R.}}
(\byear{1996}).
\btitle{Regression shrinkage and selection via the lasso}.
\bjournal{J. Roy. Statist. Soc. Ser. B}
\bvolume{58}
\bpages{267--288}.
\bid{issn={0035-9246}, mr={1379242}}
\bptok{imsref}%
\end{barticle}
%
\endbibitem

\bibitem[\protect\citeauthoryear{Tibshirani}{1997}]{T97}
%
\begin{barticle}[auto:STB|2011/11/09|09:54:39]
\bauthor{\bsnm{Tibshirani},~\bfnm{R.}\binits{R.}}
(\byear{1997}).
\btitle{The LASSO method for variable selection in the Cox model}.
\bjournal{Stat. Med.}
\bvolume{16}
\bpages{385--395}.
\bptok{imsref}%
\end{barticle}
%
\endbibitem


\bibitem[\protect\citeauthoryear{van~de Geer}{1995}]{vG95}
%
\begin{barticle}[mr]
\bauthor{\bparticle{van~de} \bsnm{Geer},~\bfnm{Sara}\binits{S.}}
(\byear{1995}).
\btitle{Exponential inequalities for martingales, with application to maximum
likelihood estimation for counting processes}.
\bjournal{Ann. Statist.}
\bvolume{23}
\bpages{1779--1801}.
\bid{doi={10.1214/aos/1176324323}, issn={0090-5364}, mr={1370307}}
\bptok{imsref}%
\end{barticle}
%
\endbibitem

\bibitem[\protect\citeauthoryear{van~de Geer and B{\"u}hlmann}{2009}]{vGB09}
%
\begin{barticle}[auto:STB|2011/11/09|09:54:39]
\bauthor{\bparticle{van~de} \bsnm{Geer},~\bfnm{S.}\binits{S.}}
\AND
\bauthor{\bsnm{B{\"u}hlmann},~\bfnm{P.}\binits{P.}}
(\byear{2009}).
\btitle{On conditions used to prove oracle results for the LASSO}.
\bjournal{Electron. J. Stat.}
\bvolume{3}
\bpages{1360--1392}.
\bptok{imsref}%
\end{barticle}
%
\endbibitem

\bibitem[\protect\citeauthoryear{Wang et~al.}{2009}]{WNZZ09}
%
\begin{barticle}[mr]
\bauthor{\bsnm{Wang},~\bfnm{S.}\binits{S.}},
\bauthor{\bsnm{Nan},~\bfnm{B.}\binits{B.}},
\bauthor{\bsnm{Zhou},~\bfnm{N.}\binits{N.}} \AND
\bauthor{\bsnm{Zhu},~\bfnm{J.}\binits{J.}}
(\byear{2009}).
\btitle{Hierarchically penalized {C}ox regression with grouped variables}.
\bjournal{Biometrika}
\bvolume{96}
\bpages{307--322}.
\bid{doi={10.1093/biomet/asp016}, issn={0006-3444}, mr={2507145}}
\bptok{imsref}%
\end{barticle}
%
\endbibitem

\bibitem[\protect\citeauthoryear{Wu and Lange}{2008}]{WL08}
%
\begin{barticle}[mr]
\bauthor{\bsnm{Wu},~\bfnm{Tong~Tong}\binits{T.~T.}} \AND
\bauthor{\bsnm{Lange},~\bfnm{Kenneth}\binits{K.}}
(\byear{2008}).
\btitle{Coordinate descent algorithms for lasso penalized regression}.
\bjournal{Ann. Appl. Stat.}
\bvolume{2}
\bpages{224--244}.
\bid{doi={10.1214/07-AOAS147}, issn={1932-6157}, mr={2415601}}
\bptok{imsref}%
\end{barticle}
%
\endbibitem

\bibitem[\protect\citeauthoryear{Yuan and Lin}{2007}]{YL07}
%
\begin{barticle}[mr]
\bauthor{\bsnm{Yuan},~\bfnm{Ming}\binits{M.}} \AND
\bauthor{\bsnm{Lin},~\bfnm{Yi}\binits{Y.}}
(\byear{2007}).
\btitle{On the non-negative garrote estimator}.
\bjournal{J. R. Stat. Soc. Ser. B Stat. Methodol.}
\bvolume{69}
\bpages{143--161}.
\bid{doi={10.1111/j.1467-9868.2007.00581.x}, issn={1369-7412}, mr={2325269}}
\bptok{imsref}%
\end{barticle}
%
\endbibitem

\bibitem[\protect\citeauthoryear{Zhang}{2010}]{Z09}
%
\begin{barticle}[mr]
\bauthor{\bsnm{Zhang},~\bfnm{Cun-Hui}\binits{C.-H.}}
(\byear{2010}).
\btitle{Nearly unbiased variable selection under minimax concave penalty}.
\bjournal{Ann. Statist.}
\bvolume{38}
\bpages{894--942}.
\bid{doi={10.1214/09-AOS729}, issn={0090-5364}, mr={2604701}}
\bptok{imsref}%
\end{barticle}
%
\endbibitem

\bibitem[\protect\citeauthoryear{Zhang and Huang}{2008}]{ZH08}
%
\begin{barticle}[mr]
\bauthor{\bsnm{Zhang},~\bfnm{Cun-Hui}\binits{C.-H.}} \AND
\bauthor{\bsnm{Huang},~\bfnm{Jian}\binits{J.}}
(\byear{2008}).
\btitle{The sparsity and bias of the {LASSO} selection in high-dimensional
linear regression}.
\bjournal{Ann. Statist.}
\bvolume{36}
\bpages{1567--1594}.
\bid{doi={10.1214/07-AOS520}, issn={0090-5364}, mr={2435448}}
\bptnote{check year}%
\bptok{imsref}%
\end{barticle}
%
\endbibitem

\bibitem[\protect\citeauthoryear{Zhao and Yu}{2006}]{ZY06}
%
\begin{barticle}[mr]
\bauthor{\bsnm{Zhao},~\bfnm{Peng}\binits{P.}} \AND
\bauthor{\bsnm{Yu},~\bfnm{Bin}\binits{B.}}
(\byear{2006}).
\btitle{On model selection consistency of {L}asso}.
\bjournal{J. Mach. Learn. Res.}
\bvolume{7}
\bpages{2541--2563}.
\bid{issn={1532-4435}, mr={2274449}}
\bptok{imsref}%
\end{barticle}
%
\endbibitem

\end{thebibliography}
\end{document}